\documentclass{article} 
\usepackage{prooftree} 
\usepackage{etex} 
\usepackage[utf8]{inputenc}
\usepackage[OT1]{fontenc}
\usepackage[francais]{babel}
\usepackage{amsmath, amssymb, amsthm}
\usepackage{amstext, amsfonts, a4}
\newcommand\seq\vdash 
\usepackage{gastex}
\usepackage{xspace}

\usepackage{epic}
\usepackage{eepic}
\usepackage{psfrag}

\newcommand\infi{\mathit{inf}}
\usepackage{amscd}
\usepackage{graphics,graphicx,subfigure}
\usepackage{comment}
\usepackage{epsfig,epsf}
\usepackage{float} 

\newtheorem{Theo}{Th\'eor\`eme}[section]

\newtheorem{Propo}{Proposition}[section]

\newtheorem{Rema}{Remarque}[section]

\newtheorem{Exem}{Exemple}[section]

\def\k#1{\kern#1em}
\def\Ib#1{{I\kern-.35em#1}}
\def\Ibb#1{{I\kern-.23em#1}}

\def\({\left(}
\def\){\right)}
\def\[{\left[}
\def\]{\right]}
\def\<{\langle}
\def\>{\rangle}

\def\ogg~{{\rm \og}}

\newcommand{\fnmark}{\mbox{$^{\arabic{footnote}}$}}
\makeatletter
\renewcommand{\@makefntext}[1]%
   {\noindent\makebox[1.8em][r]{\fnmark}#1}
\makeatother

\def\N{\mathbb{N}}

\newcommand{\fl}{\rightarrow}

\def\overset#1#2{\mathrel{\mathop{\kern0pt #2}\limits^{#1}}}

\def\k#1{\kern#1em}
\def\Ib#1{{I\kern-.35em#1}}
\def\Ibb#1{{I\kern-.23em#1}}

\usepackage{a4wide}

\usepackage{xypic}
\xyoption{all}
\input xymatrix

\usepackage{proof}
\usepackage{tikz}
\usetikzlibrary{arrows,automata,decorations.pathmorphing,backgrounds,positioning,fit,shapes,arrows}

\newcommand\ma[1]{\emph{``#1''}} 
\newcommand\ttt{\mathsf{t}} 
\newcommand\tttt{\mathsf{t_v}} 
\newcommand\eee{\mathsf{e}} 
\newcommand\vvv{\mathsf{v}} 
\newcommand\event{\mathsf{v}} 
\newcommand\imp\supset

\begin{document}


\noindent\textit{\large
L'éditeur \textbf{Cassini} publiera une version remaniée de ce texte dans un volume de \textbf{Leçons de mathématiques d'aujourd'hui} coordonné par \textbf{Géraud Sénizergues}.} 

\bigskip 
\bigskip 

\begin{center} \huge \bf 

Logique mathématique\\ et linguistique formelle\\[2ex]  

\normalfont \normalsize 
\large

\textbf{Leçon de mathématiques d'aujourd'hui du 7 juillet 2011} \\[1ex]

donnée par \textbf{Christian Retoré} (Université de Bordeaux)\\[1ex]

rédigée par \textbf{Noémie-Fleur Sandillon-Rezer} (Université de Bordeaux). 

\end{center} 

\vfill 

\begin{abstract} \normalsize 
De par son étymologie même, la logique est intimement liée au langage, comme en témoignent mathématiciens et philosophes antiques et médiévaux. Au début du XXe siècle, la crise des fondements des mathématiques a inventé la logique mathématique et l'a imposée comme un fondement des mathématiques centré sur le langage. Que sont devenus, dans ce nouveau cadre mathématique, les liens unissant logique et langage naturel? Après un aperçu de l'histoire de la connexion entre logique et linguistique traditionnellement focalisée sur la sémantique, nous ferons ensuite un gros plan sur quelques questions actuelles: 
\begin{itemize} 
\item la grammaire vue comme un système déductif 

\item le passage de la  structure grammaticale d'une phrase
à une formule logique exprimant son sens 

\item la prise en compte du contexte lors de l'interprétation 
\end{itemize} 
L'exposé montrera comment la théorie des types constitue un cadre adapté à décrire la grammaire du langage naturel et son interprétation à tout niveau (mots, phrases, discours).
\bigskip 
\end{abstract} 

\section{Introduction}
Même s'il est possible de trouver quelques précédents, on peut situer l'origine des liens entre logique et  langage dans les travaux d'Aristote (384--322 av. J.-C.), et jusqu'au Moyen-Âge au moins les deux sujets étaient totalement liés. 
Ce n'est qu'au début du XX\up{e} siècle, lors de la crise des fondements des mathématiques, que la logique mathématique telle que nous la connaissons s'est mise en place et qu'elle s'est imposée 
comme fondement des mathématiques mais aussi comme fondement de la linguistique formelle, en tout cas de la sémantique du langage naturel. C'est à ce même moment qu'on a commencé à distinguer la logique de l'étude du langage humain. 

La logique mathématique peut être vue de deux manières :
\begin{itemize} 
\item 
 l'étude de la  logique 
à l'œuvre dans les mathématiques 
dont les raisonnements constituent pour  Aristote le canon de la logique. Il s'agit alors de  \emph{logique appliquée aux mathématiques}, 
\item 
et l'étude formelle de la logique ordinaire, 
avec des structures mathématiques connues, il s'agit alors de  \emph{mathématiques appliquées à la logique}. 
\end{itemize} 
Cette leçon s'ancre plutôt dans le second courant. Avec les outils mathématiques les plus classiques possibles, si possible élémentaires, peut-on dire quelque chose sur les raisonnements, ou la vérité 
et sur leur énonciation dans le langage commun?  

\newpage 

Après un aperçu de l'histoire de la connexion entre logique et linguistique
traditionnellement focalisée sur la sémantique, c'est-à-dire sur le sens de ce que nous disons, 
nous ferons ensuite un gros plan sur quelques questions actuelles:
\begin{itemize}
\item la grammaire vue comme un système déductif
\item le passage d'une structure grammaticale à une formule logique exprimant son sens
\item la prise en compte du contexte lors de l'interprétation
\end{itemize}
La leçon montrera comment la théorie des types constitue 
un cadre adapté à décrire la grammaire du
langage naturel et son interprétation sémantique à tout niveau (mots, phrases,
discours).

\section[Une histoire des liens entre logique et langage 
où la sémantique a la part belle]
{Une histoire des liens entre logique et langage\\ où la sémantique a la part belle}

Au départ, la logique, tout comme les mathématiques, se travaillent sans langage formel, tout se fait uniquement en utilisant la rhétorique, l'art oratoire. Ce procédé naturel est encore couramment utilisé: la présentation de calculs ou de raisonnements mathématiques, en particulier à l'oral,  n'utilise des formules qu'au bon moment, généralement en les écrivant au tableau, et non à tout bout de champ. 
Du reste, à l'origine il n'y avait pas ces dites formules qui peuvent arriver à point nommé pour éclaircir un énoncé ou un raisonnement --- rappelons que les variables sont seulement apparues à la renaissance avec François  Viète alors que la pratique des mathématiques avait plus de vingt siècles! 

Même si grammaire (syntaxe) et logique (sens) sont intimement liées jusqu'au 19e siècle, les études logiques sur le langage semblent plutôt porter sur la sémantique.

\begin{description}
\item[Logos : le discours rationnel.] Ce mot dont dérive le mot ``logique" n'inclut pas l'étude de l'Odyssée, ou d'autres romans. Le modèle de discours utilisé est celui mis en place par Platon, voire Zénon, dans leurs  démonstrations mathématiques, en particulier celles de l'école de Pythagore, considérées comme le canon d'un raisonnement bien construit, ``carré". 
Le but  de la logique est de donner au raisonnement courant la rigueur mathématique présente dans une preuve pythagoricienne.
\end{description}

\subsection{D'Aristote aux penseurs médiévaux}
L'opinion commune  attribue non sans raison la fondation de la logique à Aristote, trois siècles avant Jésus-Christ.  Les bases de cette logique sont les \emph{catégories}, qui correspondent aux prédicats, et les \emph{termes} qui en fait correspondent assez mal aux termes de la logique contemporaine. Combinés entre eux, ces éléments donnent lieu à des \emph{propositions}, qui peuvent être vraies ou fausses.

\paragraph*{Termes et catégories}  
C'est une vision assez structurée d'un univers avec plusieurs sortes de propriétés. Cette tendance était assez courante, elle a été  utilisée par d'autres philosophes importants comme Porphyre mais elle s'est un peu perdue avec le temps, et il n'en reste rien ni dans la théorie des ensembles ni chez  Frege (voir plus loin paragraphe \ref{frege}) : par exemple, 
un nombre réel est un ensemble suivant les fondements ensemblistes des mathématiques --- mais qui le voit ainsi?  Les affirmations sont divisées, selon Aristote,  en une dizaine de  catégories:
\footnote{Le mots catégories peut être traduit par ``chef d'accusation", c'est un terme juridique, c'est ce dont on peut accuser un ``terme": Socrate (un terme particulier) peut être accusé d'être mortel (une catégorie).  En latin, l'expression désignant les catégories était ``praedicamentes". Notons qu'aujourd'hui, une catégorie au sens commun désigne plutôt une sorte, un ensemble de personnes, c'est-à-dire ceux pour lesquels le prédicat est vrai, interprétation déjà présente chez Porphyre. }
 \newline 

\begin{center}
\begin{tabular}{|l|l|}
\hline
\textbf{Catégories} & \textbf{Lexique, ontologie}\\
\hline
Substance (ou essence) & homme ou cheval \\
Quantité & long de deux coudées\\
Qualité & blanc, grammairien\\
Relation & double, maître, savoir\\
Lieu & au Forum\\
Temps & hier\\
Position & assis\\
Possession & armé\\
Action & couper\\
Passion & couple\\
\hline
\end{tabular}\newline
\end{center}

A cette notion, il faut ajouter les termes, qui sont soit \emph{généraux}, soit \emph{particuliers} : ``homme" est un terme général et ``Socrate" un terme particulier par exemple. Etant donné que nous sommes déformés par notre manière de raisonner en logique, nous avons tendance à voir très différemment ces deux termes, l'un  comme une variable parcourant un ensemble et l'autre comme une constante, mais, selon la logique de l'époque les deux sont uniquement des termes, l'un désignant un individu particulier, Socrate, et l'autre désignant un ensemble, voire un individu générique. 

\paragraph*{Propositions et syllogismes}
Les propositions sont des phrases particulières de quatre sortes : \newline

\begin{center}
\begin{tabular}{|ll|l|}
\hline
\textbf{Type} && \textbf{Exemple}\\
\hline 
A universelle affirmative &: & Tous les hommes sont mortels.\\
I particulière affirmative &: & Certains hommes sont philosophes. \\
E universelle négative &: & Aucun homme n'est immortel.\\
O particulière négative &: & Tous les hommes ne sont pas philosophes\footnotemark \\
\hline
\end{tabular}\newline
\end{center}
\footnotetext{Le \textit{pas tous} n'est pas lexicalisé \ma{Tous les hommes ne sont pas philosophes.}, mais il peut être formulé  différemment \ma{Certains hommes ne sont pas philosophes.}. Cependant, ces deux phrases ont des significations subtilement différentes: la première parle de l'ensemble des hommes, tandis que la seconde place l'accent sur ceux d'entre eux qui ne sont pas philosophes.}

A partir de ces quatre sortes de phrases sont définis  les célèbres syllogismes\footnote{Un syllogisme est une formalisation du raisonnement qui consiste à prendre deux propositions, ou prémisses, et en déduire une conclusion.}, qui sont tous répertoriés, classifiés avec les lettres qui donnent la forme des propositions. Nous ne donnons pas la nomenclature détaillée, mais signalons toutefois qu'elle consiste en des noms facilitant la mémorisation de la construction logique, noms dont  seules les voyelles sont prises en compte. Par exemple ``Baroco" est le nom donné au syllogisme dont les deux prémisses sont de type ``A" et ``O", soit universelle affirmative, et particulière négative, et dont la conclusion est elle aussi de type particulière et négative.

\begin{Exem} Le syllogisme Baroco 
\begin{tabular}[t]{lc}
Tout parisien est pressé & A \\ or certains conducteurs ne sont pas pressés & O\\ \hline donc certains conducteurs ne sont pas parisiens. & O
\end{tabular}
\end{Exem}

Les principes du raisonnement, les opérations possible sur les propositions,  sont généralement associés au nom Aristote, mais il est probable que certains commentateurs et traducteurs aient pris des libertés. 

En tout cas, voici deux principes usuellement attribués à Aristote : 
\begin{description}
\item[Non contradiction] : $\neg (p \wedge \neg p)$ \footnote{A des religieux qui contestaient l'utilité du principe  de non contradiction, car certains énoncés d'un dogme religieux peuvent sembler contradictoires, Avicenne répondit : ``Toute personne niant le principe de non contradiction devrait être battue et brulée jusqu'à ce qu'elle admette qu'être battue n'est pas la même chose que ne pas être battue, et qu'être brulée n'est pas la même chose que ne pas être brûlée".}
\item[Le tiers exclu] : $ p \vee \neg p$\footnote{Le tiers exclu est remis en cause par l'intuitionnisme, car il produit des preuves non constructives. Néanmoins,  dans des mondes finis, c'est un principe accepté de tous --- les semi-intuitionnistes français Baire, Borel et Lebesgue l'acceptent pour les ensembles dénombrables.}
\end{description}

\subsection{Leibniz}
Pour Leibniz, nos pensées sont la ``somme" de plus petites pensées. Selon lui, il convient de décomposer toutes les formules en sortes de propositions atomiques, et de construire des raisonnements plus complets à partir de ces propositions atomiques. La logique est alors vue comme un calcul (que l'on peut voir comme l'algèbre de Boole  ou les lois de De Morgan). Ce même Leibniz a introduit une notion d'égalité que l'on retrouve en théorie des types et dans les langages de programmation associés: 

\begin{Propo}[Egalité de Leibniz :] 
Si on peut, dans chaque proposition, remplacer $x$ par $y$ en préservant la vérité, $x = y$. Autrement dit : si deux objets ont les mêmes propriétés, ils sont égaux.
\end{Propo}
Cette définition ne correspond pas toujours à l'égalité.  Si le langage est très réduit,  par exemple  s'il n'y a que deux prédicats \emph{pair}, et \emph{impair}, rien ne permet d'exprimer la différence entre deux nombres pairs comme 8 et 10: ils sont égaux au sens de Leibniz, relativement à un langage, mais cela ne veut pas dire qu'ils soient absolument égaux, ils ne le sont qu'au regard du référentiel langagier choisi pour parler des nombres. 

Leibniz est aussi considéré comme l'inventeur supposé des Mondes Possibles. 
Cela dit, on trouve des références aux mondes possibles chez des philosophes tels que Al Farabi et Averroes  (\textit{Notre monde est le meilleur des mondes possibles}).

\subsection{Frege}\label{frege}
En tant qu'être humain, Frege n'est guère recommandable, puisque vers la fin de sa vie, en 1925,  son journal contient quelques diatribes antisémites. Mais il a  cependant apporté une certaine révolution à la logique. Il prend pour parti de tout réduire à la logique (position appelée Logicisme), en commençant par l'arithmétique. Pour cette dernière, il est le premier à montrer qu'il est possible de définir les nombres entiers dans un système logique. 
 
Frege est un grand adepte du formalisme. Cela a de bons côtés: il est le premier à écrire des formules avec les quantificateurs que nous connaissons aujourd'hui alors notés $(x)$ (aujourd'hui $\forall x$) et $Ex$ (aujourd'hui $\exists x$). Grâce à cela, 
les formules suivantes, qui étaient  indiscernables chez Boole, sont clairement distinctes chez Frege: 

\begin{description}
\item  $\forall x[I(x) \rightarrow (E(x) \vee A(x))]$
\item $\forall x (I(x) \rightarrow E(x)) \vee \forall x (I(x) \rightarrow A(x))$
\end{description} 

Cette notation permet aussi  d'écrire convenablement les ambiguïtés de portées:
\begin{itemize} 
\item 
les enfants mangent une pizza 
\item 
$\forall x\ \exists y\ \mbox{pizza}(y) \land [\mbox{enfant}(x) \rightarrow  \mbox{mange}(x,y))]$
\item 
$\exists y\ \mbox{pizza}(y) \land  \forall x\   [\mbox{enfant}(x) \rightarrow  \mbox{mange}(x,y))]$
\end{itemize}

Frege  a aussi introduit la théorie naïve des ensembles, même si  certains paradoxes y ont été trouvés, dont celui de Russell : 

 \begin{Exem}[Paradoxe de Russell] L'ensemble des ensembles qui ne s'appartiennent pas: 
  $X = \{u\mid u \notin u\}$.  A-t-on  $X \in X$ ou $X\not\in X$. 
 Si $X\in X$ alors $X\not\in X$ et si $X\not\in X$ alors 
  $X\in X$. Un tel ensemble n'existe pas. 
  \end{Exem}

Sur ce sujet, Ernst Zermelo, mathématicien allemand qui oscilla entre mathématiques et philosophie, pense qu'on peut toujours définir  un ensemble par des formules, à condition de se restreindre aux éléments d'un ensemble déjà défini, comme par exemple avec la formule : $X = \{u\in Y \mid P(u) \}$, $P$ étant un prédicat quelconque mais bien défini sur les objets de l'ensemble $Y$. Les formules logiques servent en fait à définir un sous ensemble d'un ensemble, et non à définir un ensemble en général. Le paradoxe de Russell  est ainsi évité. 

\subsection{Quelques notions logiques toujours pertinentes en linguistique}
Jusqu'à la fin du Moyen-Âge, la logique est la \textit{scientia sermocinalis}, la science du langage. C'est l'art  de dire les choses, l'organisation des mots et des phrases, et l'étude  du sens que  véhiculent les suites de mots et de phrases ainsi organisée.

\begin{description}
\item[Les Modalités] comme ``il est nécessaire que" ou ``il est possible que", initialement dues à Aristote  sont toujours d'actualité, par exemple dans la recherche sur la spécification et la vérification de  processus informatiques, en particulier parallèles. 
\item[Sens et dénotation.]  ``L'étoile du matin" et ``l'étoile du berger" sont deux termes qui dénotent Vénus, mais leur sens est différent. Le sens est, dans ce cas là, le cheminement effectué pour identifier l'objet. Ce n'est pas un cheminement psychologique, mais un cheminement dans une connaissance commune à la communauté linguistique concernée. Cette distinction est le titre d'un ouvrage de Frege, mais certains attribuent l'origine de cette distinction à Abélard. 
\item[Ambiguités :] Il y a toutes sortes d'ambiguïtés: \\ 
\begin{tabular}{|p{2cm}|p{6.3cm}|p{6.3cm}|}
\hline
\textbf{Sens} & J'ai fini mon livre. & Fini de le lire ou de l'écrire ? \\
& C'est un avocat marron. & fruit ou membre du barreau? \\
\textbf{Structure} &Il regarde quelqu'un avec des lunettes noires. & Qui porte les lunettes noires ? \\
\textbf{Portée} & Les enfants prendront une pizza. & Une pizza par enfant ou une seule pour tous les enfants ? \\
\textbf{Implicite} & Pierre aime sa femme; Jean aussi. & La femme de qui Jean aime-t-il ? \\
\hline
\end{tabular}
\item[Prédicats vagues :] Ils sont très durs à modéliser dans un cadre. Par exemple, où commence la couleur bleue ? Le turquoise est il encore bleu ou déjà vert ? A partir de quelle taille quelqu'un est il grand? 
Cela conduit au paradoxe des Sorites. Si $X$ est grand et qu'on lui retire $1$ le résultat est grand mais ce principe conduit à trouver grand n'importe quel nombre, y compris $0$! 
\item[Mots catégorématiques :]  Ce sont les termes comme définis par Aristote, cf. ci-dessus. 
\item[Mots syntacatégorématique :] Ce sont les mots qui ont une fonction logique (\emph{tous, aucun, si, non,et,...} )
\item[Suppositio (d'un terme) :] principes médiévaux régissant le sens en contexte : le mot va être sujet d'un prédicat. On différencie le sens propre (au nombre de deux) et le sens impropre (figures de style telles que les métaphores, les métonymies...).
\item
\begin{tabular}{|p{2cm}|p{6.45cm}|p{6.45cm}|}
\hline
\textbf{Sens propre}\footnotemark & Vélo a quatre lettres. & On parle du mot ``vélo" (suppositio materialis) \\
& Un vélo lui a été offert. & On parle de l'objet ``vélo" (suppositio formalis)\\
\hline
\textbf{Sens impropre} & Un vélo est pratique dans Bordeaux. & Le moyen de transport est pratique, pas l'objet.\\
& Son vélo est crevé. & Ce n'est pas le vélo qui est crevé mais une de ses roues. \\
\hline
\end{tabular}
\footnotetext{On rappellera qu'en latin, il n'y a pas d'article.} 

\item[Querelle des universaux :] Cette querelle porte sur le statut des termes : un terme général n'est-il qu'un nom, une sorte d'abréviation pour tous les individus qui tombent sous ce concept (nominalisme) ou a-t-il une existence propre (réalisme) ? Cette question est reliée à celle des éléments génériques, lesquels  capturent l'essence d'une classe. Dans un contexte mathématique, en disant soit ``$n$ un entier" pour prouver une propriété sur les entiers, le $n$ utilisé dans la preuve existe-t-il vraiment ? C'est un entier générique, il n'existe pas vraiment en tant que ``$n \in \mathbb{N}$", mais son usage est compréhensible voire naturel pour un mathématicien (mais peut-être ne l'est-il pas pour une partie non négligeable de nos étudiants). 

\item[Lecture \emph{de dicto} et \emph{de re}] Cette question introduite par Thomas d'Aquin et ainsi nommée par Guillaume d'Ockham est  en fait un problème de portée entre un quantificateur située dans une proposition subordonnée et le même quantificateur portant sur la proposition principale: \ma{James Bond croit qu'un chercheur de l'IMB est un espion} a deux lectures possibles, et la portée du ``il existe" diffère dans les deux cas :
\item
\begin{tabular}{|p{2cm}|p{6.45cm}|p{6.45cm}|}
\hline
\textbf{De re} &  James Bond a vu Gilles effectuer une action louche, et le soupçonne & il existe quelqu'un dont James Bond pense qu'il est un espion\\
\hline
\textbf{De dicto} & James Bond a trouvé un indice qui laisse supposer que quelqu'un est un espion &James Bond pense qu'il existe un espion. \\
\hline
\end{tabular}
\end{description}

Afin de personnaliser les notions introduites au début du XXe siècle, on pourrait relier les notions passées en revue à certains logiciens: 

\begin{description}
\item[Formules :] Frege;
\item[Preuves formelles :] Frege, Hilbert, Herbrand, Gentzen, G\"odel; 
\item[Modèles :] Löwenheim, Skolem, Tarski, Herbrand, G\"odel; 
\item[Théorie des ensembles :] Zermelo, Fränkel, Gödel;\footnote{Et oui, G\"odel apparaît un peu partout. 
Sa thèse a établi la complétude du calcul des prédicats (le théorème d'Herbrand peut aussi se lire ainsi, lorsque le langage est dénombrable): tout ensemble d'énoncé non contradictoire admet un modèle. Il a aussi étudié les modèles de la théorie des ensembles en montrant que l'axiome du choix et l'hypothèse du continu sont compatibles avec la théorie des ensembles. Et il a montré qu'il existe des propriétés vraies des entiers qui ne sont pas démontrables à partir des axiomes de l'arithmétique (incomplétude).}.
\end{description}

\subsection{Une vision logique du sens est-elle possible} 

 Les structures logiques évoquées durant cette rapide histoire
 permettent-elles de modéliser la signification? 
Plusieurs visions du sens d'une phrase sont possibles, et sans aller dans les détails donnons-en quelques unes. 

Les \emph{empiristes}, par exemple, vont l'approcher du point de vue des ``contenus mentaux suscités par des signes" (Locke, Berkeley, Hume), ce qui est quasiment impossible à formaliser. \\
Il y a aussi l'approche axée sur l'\emph{usage et l'interaction} (Wittgenstein, Searle, Brandom) : le sens d'une phrase est composé de toutes les choses qui y répondent. Du point de vue d'un dialogue, ce sont toutes les bonnes réponses, la continuation du dialogue. \\
La méthode de \emph{vérification/réfutation} (Quine) se définit sur la possibilité de vérifier ou de réfuter une phrase. \\
Le \emph{pragmatisme}voit le sens comme l'effet du dire sur l'environnement (Peirce). \\

Dans cet exposé, faute de temps, les seules interprétations seront formulées dans des modèles en termes de   \emph{conditions de vérité} (Frege, Tarski, Davidson) et c'est assurément regrettable car d'autres interprétations sont plus fines et plus pertinentes pour l'étude du langage. 
Les conditions de vérités permettent de dire
quand une formule est vraie, et quand elle ne l'est pas notamment par une \emph{ l'interprétation dans des monde(s) possibles} (Putnam, Kripke ont mathématisé cette notion).

\section{Formules, preuves et modèles} 
 \subsection{Logique propositionnelle}
 
 \subsubsection{Formules}  
Commençons pas  déterminer ce qu'est la logique propositionnelle en tant que langage formel:  
les formules sont fabriquées à partir de propositions atomiques en utilisant des connecteurs logiques:  
\begin{description}
\item[Propositions atomiques :] $p_i$, $i \in \N$
\item[Connecteurs usuels :] $A \vee B$, $A \wedge B$, $A \rightarrow B$, $A \& B$, ...
\end{description}

\subsubsection{Preuves} 
Et avec la notion de formule va celle de preuve, sous forme d'arbre, comme montré sur la figure \ref{proof}. L'arbre de preuve a pour racine la formule à prouver, les noeuds ou branchement représentent les étapes des preuves appelées règles. 

Les feuilles  sont appelées hypothèses, à moins qu'elles ne soient annulées (on dit aussi déchargées) au cours de la règle $\Rightarrow_{i}$ qui introduit le connecteur $\Rightarrow$: selon cette règle, si  $B$ a été établi en utilisant $n$ fois l'hypothèse $A$ alors  $A\Rightarrow B$ a été établi 
sans utiliser l'hypothèse $A$, et les occurrences correspondantes de $A$ sont mises entre crochets en indiquant la règle qui a causé leur annulation  --- la règle dit en fait que dans ces circonstances $A\Rightarrow B$ a été établi 
 en n'utilisant qu'un sous ensemble quelconque des occurrences de l'hypothèse $A$, ce qui n'est pas faux; cependant, dans la pratique usuelle du raisonnement, toutes les occurrences de l'hypothèses $A$ sont déchargées.\footnote{Dans les présentations à la Hilbert des preuves, en particulier dans les écrits de Herbrand, ce que nous appelons ici la règle $\Rightarrow_{i}$  est le  \emph{théorème de la déduction}, relativement difficile à établir.}
 
 La règle duale est connue de tous: il s'agit du \emph{modus ponens}. Si $A$ a été établi avec certaines hypothèses $\Gamma$ et $A\Rightarrow B$ avec certaines hypothèses $\Delta$ alors  $B$ peut être établi sous les hypothèses $\Gamma$ et $\Delta$. 
 
\begin{figure}[t] 
\label{proof} 
\begin{prooftree} 
	\[ \begin{array}{c}
	   \ldots \ldots [A]_\ell \ldots [A]_\ell \ldots A 
	   \end{array} 
	\leadsto 
	B 
	\using 
	\]
	\justifies 
	A \Rightarrow B
	\using \mbox{$\Rightarrow_{i} \ell$\ certains $A$ ne sont plus des hypothèses} 
\end{prooftree} 
\hspace{2cm} 
\begin{prooftree} 
	\[\Gamma
	\leadsto 
	A\Rightarrow B
	\] 
	\[\Delta 
	\leadsto 
	A
	\] 
	\justifies 
	B
	\using\Rightarrow_{e}
\end{prooftree} 
\caption{Règles d'introduction et d'élimination de l'implication} 
\end{figure}

Qu'est ce qu'une preuve de $F$? c'est une preuve sans hypothèse, c'est-à-dire une preuve dont toutes les hypothèses ont été annulées au cours de la preuve. La plus simple est assurément celle que voici: 

\begin{center} 
\begin{prooftree}
[A]_1
\justifies 
A\Rightarrow A
\using \Rightarrow_i 1
\end{prooftree} 
\end{center} 

Il est possible d'éviter les arbres et les liens (le ``$\ell$" de la règle $\Rightarrow_i$ vers les hypothèses, il faut alors écrire à chaque étape de la preuve toutes les hypothèses non annulées à cet instant, comme montré sur figure \ref{F2} gnomettant les termes écrits en minuscules. La preuve ne peut être fidèlement représentée par des des termes simples, il y a  besoin d'une opération liant les hypothèses, la lambda abstraction: cette notion nous sera fort utile par la suite. 

La même hypothèse peut être utilisée sans avoir été prouvée, notre preuve sera  complète modulo la justification de cette hypothèse. Les règles permettent ce processus assez naturel dans le raisonnement mathématique.  On utilise $A$ \emph{ad libitum} pour établir $B$ et on en conclut $A\Rightarrow B$ \emph{sans hypothèse}. On démontre à part $A$, également sans hypothèse. Et on en conclut $B$, toujours sans hypothèse, par \emph{modus ponens}.  Cela fonctionne  comme un lemme dans une preuve mathématique usuelle. Lorsque l'on a la preuve du lemme, il suffit de l'insérer là où on a utilisé le lemme pour que notre preuve soit toujours complète et élémentaire (en ce qu'il n'y a plus de lemme), ce qui en veut pas dire qu'elle soit plus simple à suivre ou plus courte, bien au contraire.


\begin{figure}[b] 
\begin{center} 
\begin{prooftree}
\[
\vec x{:}\Gamma, z{:}A, z{:}A, x{:}A \seq t{:}B 
\justifies
\vec x{:}\Gamma, z{:}A, x{:}A \seq t{:}B 
\using contr. 
\]
\justifies 
\vec x{:}\Gamma\seq (\lambda z{:}A.\ t){:}A\Rightarrow B\qquad 
\using \Rightarrow_i
\end{prooftree} 
\qquad 
\begin{prooftree} 
\vec x{:}\Gamma\seq f{:}A\Rightarrow B\qquad 
\vec y{:}\Delta\seq u{:}A
\justifies
f(a){:}B
\using 
\Rightarrow_e
\end{prooftree} 
\end{center} 
\caption{Les même règles, introduction et élimination de l'implication, avec hypothèses explicites et termes} 
\label{F2} 
\end{figure}

A partir de là, pour avoir la logique usuelle, il faut ajouter comme axiome $A\vee \neg  A$ pour toute formule $A$, ce qui correspond au tiers exclu dont nous avons parlé précédemment. C'est étonnant qu'il faille ajouter ce principe. Du reste, en le laissant de côté, on obtient une logique fort intéressante, la logique intuitionniste, qui est constructive, c'est-à-dire qui permet d'extraire des preuves existentielles des algorithmes calculant ce dont la preuve prouve l'existence. 

\subsubsection{Modèles}
 Les modèles sont une autre façon de parler de la vérité: ils permettent de vérifier, dans une situation donnée,  si une formule est vraie ou fausse, alors que les preuves permettent de dériver des formules 
 vraies en toutes circonstances. 
 
 Pour que l'interprétation soit valide, on dit aussi bien fondée, il déjà faut qu'une formule démontrable soit vraie quelles que soient les contingences, et heureusement c'est le cas. 
 Mais il y a aussi une réciproque, appelée complétude: 
 
\begin{Theo}[Théorème de complétude (Bernays, 1926)]\label{completudeprop}
Une formule vraie dans tout modèle est démontrable : la déduction capte exactement ce 
qui peut s'établir par les valeurs de vérité.
\end{Theo}

Autrement dit : si une formule est tout le temps vraie, partout, quelles que soit les circonstances, alors il en existe une preuve --- et sinon il n'en existe pas. 

\subsubsection{Logique propositionnelle modale}

Nous disons  souvent des phrases comme ``il est possible que Claire vienne", 
et pour en parler 
des opérateurs unaires appelées modalités, sont ajoutés à la logique. Les plus courants sont  ``il est possible que", et ``il est nécessaire que", représentés respectivement par  $\diamond$ et $\square$. La notation $\square p$ correspondra donc à ``il est nécessaire que p".

Des axiomes gèrent les modalités: 
\begin{itemize}
\item $\diamond p \equiv  \neg  \square \neg  p$
\item $\square p \equiv  \neg  \diamond \neg  p$
\item $\square p \rightarrow  p \square $
\item $p \rightarrow  \square \diamond p$
\end{itemize}

Pour créer un modèle, il ne faut pas juste une attribution de valeur de vérité, mais une famille, avec une relation d'accessibilité. Donnons nous des modèles: M1, M2, M3,...  qui constituent de vrais mondes, avec des choses vraies et fausses, ainsi qu'une relation d'accessibilité entre ces mondes qui  nous amène 
d'un monde à  un autre.  En fonction de la logique modale que l'on considère, la relation d'accessibilité doit avoir plus ou moins de propriétés.\footnote{On peut aussi interpréter la logique intuitionniste définie par les règles de déduction sans le tiers exclu, en demandant à ce que l'accessibilité soit transitive et réflexive.} 

\begin{Exem} Voici l'interprétation de deux phrases modales.
\ma{Il est nécessaire que Pierre vienne} sera vraie dans un monde $M$ si dans  tous les mondes accessibles à partir de $M$ "Pierre vient."

\ma{Il est possible que Pierre vienne} sera vraie dans un monde $M$ si dans l'un des mondes accessibles à partir de $M$ \ma{Pierre vient/viendra}. 
\end{Exem}

Pour cette extension, avec les modèles évoqués,  le théorème de complétude vaut toujours: les propositions vraies dans tous les modèles sont exactement les propositions démontrables.

\subsection{Calcul des prédicats}
Revenons à la logique de Frege, à son invention, le calcul des prédicats.

\subsubsection{Formules} 
Les prédicats de base sont des relations $n$-aires, en fonction du nombre d'arguments qu'ils prennent. Même si  dans le langage courant, on a tendance à typer les choses, pour Frege, il n'y a ni type ni sorte, ni catégorie : c'est le même univers pour tout le monde. C'est ainsi que la logique en vient à considérer un nombre réel comme un ensemble, ce qui est possible et cohérent mais bien peu naturel. 

\begin{Exem}
Dans le langage courant : 
\begin{itemize}
\item dormir est unaire (quelqu'un dort) ;
\item regarder est binaire (quelqu'un regarde quelque chose) ;
\item donner est ternaire (quelqu'un donne quelque chose à quelqu'un).
\end{itemize}
En mathématiques :
\begin{itemize}
\item premier(x) est unaire (x est premier);
\item plus(x, y, z) est ternaire (x + y = z).
\end{itemize}
\end{Exem}

Afin d'éviter de compliquer cet exposé, nous  laissons délibérément de côté les notions de \emph{fonction} et de \emph{terme}:
du reste il est possible d'exprimer, dans ce genre de langage, qu'un prédicat est une fonction. 
Comme par exemple pour le $plus$ mentionné ci-dessus: $\forall x\forall y\exists z plus(x,y,z)$ (toute paire d'entier a une somme) et $\forall x\forall y\forall z\forall z'plus(x,y,z)\wedge plus(x,y,z') \Rightarrow z = z'$
(deux entiers ont au plus une somme). 

\subsubsection{Preuves} 

Les déductions dans le calcul de prédicats, requièrent, outre les règles de la logique propositionnelle vues précédemment, des règles pour gérer la quantification. 

\begin{Exem}
\begin{description}
\item[]
\end{description}
\begin{tabular}{cc}
\begin{tabular}{c}
$A(x)$ \\
\hline
$\forall x A(x)$
\end{tabular}
 & $\forall_i^-$ si pas de x dans une hypothèse libre \\
 & \\
\begin{tabular}{c}
$\forall x A(x)$ \\
\hline
$A(x)$
\end{tabular}
 & $\forall_e^-$ A(u) est une variable ou une constante
\end{tabular}

Si $A(x)$  vaut sans supposer quoi que ce soit sur ce pauvre $x$,  la preuve de $A(x)$ avec cet $x$ générique va donc  fonctionner pour tout $x$.

En revanche, on peut instancier, c'est à dire déduire $A(u)$ de $\forall x A(x)$ pour  toute constante ou variable $t$ --- et même pour tout  terme $t$ si  nous avions des termes. 

Les règles pour le $\exists$ sont déductibles par dualité.
\end{Exem}

\subsubsection{Modèles}

A proprement parler, un modèle du calcul des prédicats  a pour base un ensemble $M$ (dans lequel les variables varient !) qui représente les valeurs que peuvent prendre les variables. Cet ensemble est unique quels que soient les objets traités, il n'y a qu'une seule sorte d'objets. 

Les prédicats $n$-aires sont alors des parties de $M^n$. 

\begin{Exem}
L'addition, $plus$, comme vue précédemment, va être une partie de $M^3$.

Ne sont pris en compte que les triplets de nombres pour lesquels la somme des deux premiers donne le troisième.
\end{Exem}

Il faut fixer la valeur des variables libres pour qu'une formule soit vraie ou fausse (assignation). 
Pour interpréter les formules closes, on fait varier l'interprétation des variables. 

\begin{Exem} Exemple sur les formules closes :

$\forall xP(x)$ est vrai si pour toute valeur assignée à la variable $x$ la formule $P(x)$ est vraie.

$\exists xP(x)$ est vrai sil existe une valeur qui, assignée à $x$, rende la formule $P(x)$ vraie.
\end{Exem}

Bien évidemment, comme pour le calcul des propositions, les formules démontrées sont vraies dans tout modèle, 
et les formules vraies sous certaines hypothèses sont vraies dans tout modèle où lesdites hypothèses sont vraies. Ce résultat, aisé à obtenir par induction, s'appelle la validité (\emph{soundness}) du calcul des prédicats. 
Comme pour le calcul propositionnel, il admet une réciproque, le théorème de complétude \ref{completudepred} : 

\begin{Theo}[Théorème de complétude (G\"odel, Herbrand, 1930)] 
Une formule vraie dans tout modèle est démontrable. 
\label{completudepred}
\end{Theo}

Ce résultat affirme que si une formule est tout le temps vraie, partout, quelles que soit les circonstances, il en existe une preuve.

Autrement dit,  et compte tenu de la validité définie ci-dessus, il y a coïncidence entre la notion de vérité portée par la déduction et celle portée par les interprétations en termes de valeurs de vérité.
Une formule du langage des groupes, conséquence des axiomes de groupes, sera vraie dans tous les groupes, et, réciproquement, une formule vraie dans tous les groupes est démontrable à partir des axiomes de groupes.

\subsubsection{Modalités et modèles de Kripke}

Il y a des modèles, dits de de Kripke pour le calcul des prédicats modal. 
Ces mêmes techniques développées pour les modalités standard permettent d'interpréter les verbes de croyance comme ``croire". 
Pour changer des modalités ``il est possible " et ``il est nécessaire que" disons-en deux mots: 

\begin{Exem}
Analysons la phrase \ma{Sébastien croit que Chomsky est un informaticien.} 

Cela ne signifie pas du tout que Chomsky soit informaticien, mais uniquement que Sébastien pense que Chomsky est informaticien, car Sébastien peut se tromper. 
Cette phrase est vraie quand dans tous les mondes où les croyances de Sébastien sont des faits avérés, Chomsky est un informaticien.
\end{Exem}

Un cas particulier des modèles de Kripke est  la logique intuitionniste, qui est la logique de nos règles, si nous ne leur ajoutons pas l'axiome du tiers exclus. Dans cette logique, comme 
dans la langue naturelle, la double négation est  un peu plus faible que l'affirmation. (\ma{Vous n'ignorez pas les lois} me semble moins fort que \ma{Vous connaissez les lois}). On a donc : $\neg \neg p \not \rightarrow p$ (mais $p \rightarrow \neg \neg p$). 

Des modèles de Kripke dans lesquels chaque monde possible est un modèle classique correspondent à la logique intuitionniste, et en travaillant dessus,  des structures mathématiques bien connues apparaissent: les (pré)faisceaux. 
De nouveau on a un théorème de complétude, applicable à ce cas précis de la logique intuitionniste.
Une preuve directe, simple et récente a été obtenue par un doctorant bordelais, Ivano Ciardelli. 

\section{Des phrases aux formules --- dans la logique !}

Pour l'instant, nous avons vu comment établir et interpréter des formules, et il est communément admis qu'elles correspondent à des phrases. Pourquoi pas, mais comment leur correspondent-elles? 
En fait nous avons suivi  le point de vue historique: les logiciens ont d'abord réussi à manipuler des formules,
ou plutôt des phrases dont les formules sous-jacentes sont évidentes,  avant de se poser la question du lien entre phrases et formules. Nous allons expliciter ce lien en montrant 
comment transformer une phrase en une formule --- ou en plusieurs  formules, une phrase pouvant être ambigüe. 
Pour arriver à cette formule à partir d'une phrase, il faut d'abord l'analyser grammaticalement puis appliquer des méthodes pour aller des phrases aux formules logiques. Nous allons ici montrer quelles sont ces méthodes.

Pour ce faire, nous commencerons par définir  l'analyse grammaticale / syntaxique d'une phrase, ce qui correspond à la preuve qu'elle est bien formée. La structure syntaxique\footnote{C'est l'équivalent des arbres, que l'on dessinait  avec des stylos bille de diverses couleurs, au début du collège: on séparait les groupes de mots suivant leur fonction dans la phrase, et on traçait les liens les unissant, mais cela peut se faire automatiquement, s'informatiser.}  nous donnera un arbre de preuve, dans une logique très particulière, très restreinte, manipulant des catégories qui seront à la fois grammaticales et sémantiques. 

\subsection{Calcul de Lambek}
Le système logique que nous allons utiliser pour reconnaitre et analyser les phrases est le calcul de Lambek.

Joachim Lambek est  algébriste (connu pour son livre sur les modules), catégoricien (allègre, logique, topos) mais aussi  féru de linguistique. Il a transformé les grammaires catégorielles d'Ajdukiewicz et  Bar-Hillel en un système logique très élégant. 

\subsubsection{Principes de base}
Lambek définit des catégories (ou types) de base :
\begin{itemize}
\item $S$ pour une phrase (sentence);
\item $n$ pour les noms communs (noun);
\item $np$ pour les syntagmes nominaux (noun phrase).
\end{itemize}

A partir de ces catégories de base, construisons  un ensemble de catégories :
\begin{center}
$P ::= S, n, np, ... $\\
$L_p ::= P \mid (L_p\backslash L_p )\mid (L_p /  L_p )$
\end{center}

C'est-à-dire : $P$ est l'ensemble des catégories atomiques et $L_p$ celui des catégories composées, formées au moyen des symboles ``a sur b" ($a/b$) et ``a sous b" ($a \backslash b$). Un lexique assigne alors des catégories aux mots  puis aux suites de mots grâce aux règles ci-après. Notons dès à présent qu'un mot ou une suite de mots peuvent avoir plusieurs catégories.

Dans le sens usuel, dû à Ajdukiewicz, Bar-Hillel, si un mot $u$ (ou une suite de mots) 
a la catégorie $A$ et un mot (ou une suite de mots) $v$ a la catégorie  $A\backslash B$\footnote{lire \textit{A sous B}, alors que $A/B$ se lit \textit{A sur B}}, alors la suite de mots $uv$  a pour catégorie $B$.

La réciprocité de ces principes a été introduite par Lambek. Si une suite de mots $uv$ a pour catégorie $B$, et que le mot ou la suite de mots $u$ a pour catégorie $A$, alors une catégorie du mot ou de la suite de mots $v$ est $A\backslash B$. 
Autrement dit  $v$ est de catégorie $A\backslash B$ s'il  attend un $A$ à sa gauche.

Les ``calculs" s'effectuent ensuite comme montré dans la figure \ref{lamb}.
\begin{figure}[htb]
\begin{center}
\begin{tabular}{ccc} 
u & \multicolumn{2}{c}{v} \\
\multicolumn{1}{c}{$\overline{u : A}$} & \multicolumn{2}{c}{$\overline{v :A\backslash B}$} \\
\hline
 \multicolumn{3}{c}{ uv : B} \\
&&\\
u & \multicolumn{2}{c}{v} \\
\multicolumn{1}{c}{$\overline{u : B/ A}$} & \multicolumn{2}{c}{$\overline{v : A}$} \\
 \hline
 \multicolumn{3}{c}{ uv : B}
 \end{tabular} 
\end{center}
\caption{Avec les règles AB (Ajdukiewicz Bar Hillel) les catégories peuvent se voir comme des fractions non commutatives, qui se simplifient uniquement si le A est du bon côté.}\label{lamb}
\end{figure}

\subsubsection{Règles}
\begin{figure}[htb]
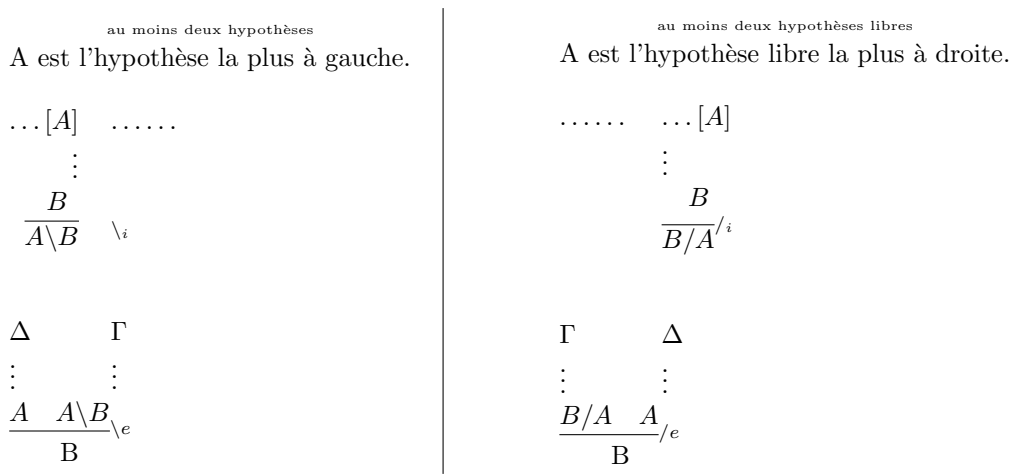

\begin{center}
\begin{tabular}{c|c}
\begin{tabular}{rl}
\multicolumn{2}{c}{\tiny{au moins deux hypothèses}} \\
\multicolumn{2}{c}{A est l'hypothèse la plus à gauche. } \\
& \\
$\dots \[A\]$ & $\dots \dots$\\
$\vdots$ & \\
\multicolumn{1}{c}{$ \quad B$}& \\
$\overline{A \backslash B}$ & $^{\backslash_i}$ \\
& \\
& \\
\multicolumn{1}{l}{$\Delta$} & $\Gamma$\\
\multicolumn{1}{l}{$\vdots$}&$\vdots$ \\
\multicolumn{2}{l}{$\underline{A \quad A \backslash B}_{\backslash e}$} \\
B 
\end{tabular} & 

\begin{tabular}{crl}
$\qquad$&\multicolumn{2}{c}{\tiny{au moins deux hypothèses libres}} \\
$\qquad$&\multicolumn{2}{c}{A est l'hypothèse libre la plus à droite.} \\
$\qquad$&& \\
$\qquad$&$\dots \dots$ & $\dots \[A\]$\\
$\qquad$&& $\vdots$\\
$\qquad$&& $\quad B$ \\ 
$\qquad$&& $\overline{B/ A}^{/_i}$ \\
$\qquad$&& \\
$\qquad$&& \\
$\qquad$&\multicolumn{1}{l}{$\Gamma$} & $\Delta$\\
$\qquad$&\multicolumn{1}{l}{$\vdots$}&$\vdots$ \\
$\qquad$&\multicolumn{2}{l}{$\underline{B/A \quad A}_{/ e}$} \\
$\qquad$&B 
\end{tabular}
\end{tabular}
\caption{On retrouve les règles AB vues auparavant (figure \ref{lamb}). 
Il y a en plus les règles d'introduction dues à Lambek. Sachant que $uv$ est de type $B$ et que $u$ est de type $A$,  $v$ est  de type $A \backslash B$.
Evidemment  les règles $/_e$ et $/_i$ procèdent symétriquement.}\label{F4}
\end{center}
\end{figure}

Les règles du calcul de Lambek sont fondées sur les grammaires AB. Elles sont résumées dans la figure \ref{F4}. Les règles $\backslash _i$ et $/_i$ sont des règles d'introduction et les règles $\backslash _e$ et $/_e$ sont des règles d'élimination.

Une autre manière de présenter les règles, avec une méthode plus ``en ligne", donne quelque chose comme  ce que nous voyons en figure \ref{d39}

\begin{figure}[htb]
	\begin{center}
		\begin{tabular}{ccc}
		$\underline{\Gamma \vdash A  \quad \Delta \vdash A\backslash B}_{\backslash e}$ & $\underline{A, \Gamma \vdash C}_{\backslash i}$ & $ \Gamma \neq \varnothing$ \\
		$\Gamma, \Delta \vdash  B$& $ \Gamma \vdash A \backslash C$& \\
		&& \\
		$\underline{\Delta \vdash B / A \quad \Gamma \vdash A}_{/e}$ & $\underline{\Gamma, A \vdash C}_{/ i}$ & $ \Gamma \neq \varnothing$ \\
		$\Delta, \Gamma \vdash B$&$ \Gamma \vdash C/A$&$$ \\
		&& \\
		\multicolumn{2}{c}{$\quad \qquad \qquad_{axiome}$}& \\
		\multicolumn{2}{c}{$\overline{A \vdash A}$}&
		\end{tabular}
	\end{center}
 \caption{Suivant  la règle $/e$,  lorsque quelque chose de catégorie $B/A$ est suivi de quelque chose de catégorie $A$, la suite des deux est de catégorie $B$.}\label{d39}
\end{figure}

On peut faire le lien avec le groupe libre : 
\begin{itemize}
\item $S$, $n$, $np$, sont les éléments de base. 
\item $(y/x)\backslash y =  (yx^{-1})^{-1}y = xy^{-1}y = x$. 
\end{itemize}
Cela rend $A \vdash B$ symétrique, ce qui n'est syntaxiquement pas vrai, le groupe libre  n'est pas un modèle complet. L'avantage de ce lien avec le groupe libre est de vérifier rapidement par un calcul que  quelque chose n'est \emph{pas} démontrable.  Si le calcul des deux membres dans le groupe libre ne donne pas le même résultat, ni $A\seq B$ ni $B\seq A$ n'est démontrable. 

\subsubsection{Une très jolie logique}
Le calcul de Lambek est une logique élégante 
où se retrouvent plusieurs caractéristiques importantes :
\begin{itemize}
\item Complétude pour des modèles algébriques ordonnés simples (une catégorie de base est interprétée par les suites de mots qui sont de cette catégorie, et cette interprétation est étendue aux catégories complexes, l'inclusion d'une catégorie $A$ dans une catégorie $B$ correspondant à la prouvabilité de $A\seq B$). (Bukowski 1982)
\item Notion de preuve normale, avec propriété de la sous formule, ce qui facilite la recherche de preuve et donc l'analyse grammaticale.  (Lambek 1958) Ladite normalisation est confluente, et s'effectue en temps linéaire. (Lambek 1958)
\item  La prouvabilité est décidable (Lambek 1958) mais NP-complète en fonction de la taille des formules à prouver (Pentus 2003)
\item Le calcul de Lambek est une sous logique de la logique intuitionniste, c'est à dire que les preuves utilisent une partie des règles de la logique intuitionniste. Cette propriété nous resservira par la suite, pour construire l'interprétation sémantique. 
\end{itemize}

\subsubsection{Grammaire de Lambek : description d'un langage ou d'une langue}
Comment fonctionne une grammaire de Lambek ?

Les règles données précédemment et présentées figure \ref{d39} et \ref{F4} sont la clé de l'édifice. Toutes les informations sur le langage  sont \emph{dans le lexique} ($lex$). Pour chaque mot, celui-ci contient toutes les catégories associées à ce mot\footnote{Un mot a parfois besoin de plusieurs catégories, comme \ma{manger} qui peut être utilisé avec ou sans complément d'objet \ma{Jean mange} ou \ma{Jean mange quelque chose}.} Quel que soit le langage $(L(lex))$, seul le lexique varie, tandis que les règles logiques restent immuables. 

Soit une phrase, vue comme une suite de mots $m_1, m_2, \dots, m_k$.  Si pour chaque mot il est possible de trouver une catégorie parmi celles que lui attribue le lexique de sorte qu'en combinant tous les types il existe une preuve de la catégorie de base $S$, la phrase est dite correcte, et appartient à $L(lex)$. 

\begin{Exem}[Exemple Simple] 
Considérons le lexique suivant, en italien pour des questions de simplification (une question se pose, en italien, sans inversion de sujet) : \\
\begin{center}
\begin{tabular}{r|ll}
Mot & Type(s) & Traduction \\
\hline
cosa & $(S/(S/np))$ & quoi \\
guarda & $(S /\infi  	)$ & il/elle regarde \\
passare & $(\infi  /np)$ & passer \\
il & $(np/n)$ & le \\
treno & $n$ & train
\end{tabular}
\end{center}

La phrase : \ma{guarda passare il treno} (il/elle regarde passer le train) donne le calcul suivant : 
\begin{center}

$(S/\infi  ) . (\infi  /np) . np/n n.$ 
\end{center}

Le type $\infi  $ représente un infinitif. Pour la simplification, regardons la figure \ref{sim}. 

\begin{figure}[htb]
	\begin{center}
$$
\infer[/e]{np}{guarda : S/\infi & \infer[/e]{\infi  }{passare:\infi  /np & \infer[/e]{np}{il:np/n & treno:n}}} 
$$
	\end{center}
 \caption{Simplification de la phrase \ma{guarda passare il treno}.  Le fait d'avoir comme catégorie racine $S$ nous assure que la phrase est bien dans $L(lex)$.}\label{sim}
\end{figure}

Comment analyser \ma{cosa guarda passare} (que regarde-t-il/elle passer ?) ?  Il  manque un groupe nominal pour pouvoir simplifier \ma{guarda passare} ($(S/\infi )\ (\infi  /np)$)  C'est là qu'interviennent les règles d'introduction ajoutées par Lambek aux grammaires AB: puisqu'avec un $np$ au bout \ma{guarda passare + $np$} est une phrase,  le groupe \ma{guarda passare} a le type $S/np$.
Pour que \ma{cosa guarda passare} soit une phrase correcte,  il faut bien attribuer à \ma{cosa} la catégorie $(S/(S/np))$.
La même preuve en remplaçant respectivement $S,\infi  ,np$ 
par $a,b,c$ établit 
 la transitivité du symbole ``$/$" dans  le calcul de Lambek : $a/b, b/c \vdash a/c$.
\end{Exem}

\begin{Exem}[Mots vides]
Une particularité du calcul de Lambek en tant que système logique est d'imposer à toute preuve ou partie de preuve d'avoir toujours au moins une hypothèse : on ne peut considérer de mot vide (alors que selon une logique usuelle, un mot vide devrait avoir  pour catégorie toute formule logique vraie sans hypothèse). Il faut qu'il en soit ainsi, car sans cela,  il serait possible de considérer \ma{un très exemple} comme un groupe nominal (voir figure \ref{void}). 
\begin{figure}[htb]
	\begin{center}
$$
\infer[/e]{np}{un : np/n & \infer[/e]{n}{\infer[/e]{n/n}{tres : (n/n)/(n/n) & \infer[/i-\alpha]{n/n}{[n]_\alpha} }  & exemple : n} } 
$$
	\end{center}
 \caption{Les adjectifs prennent en argument un nom commun et les retransforment en une sorte de nom commun (\ma{petite maison} et \ma{maison} vont se traiter, d'un point de vue syntaxique, de la même manière). Les adjectifs antéposés ($n\backslash n$) sont donc des mots qui prennent à leur droit un nom ($n$) et qui rendent un nom ($n$). En acceptant un mot vide en guise d'adjectif, \ma{un très exemple} est bien  de type $np$, qui va pouvoir s'insérer dans une phrase, alors que ce n'est clairement pas un groupe nominal valide. C'est pour cela qu'il doit toujours y avoir une hypothèse non vide.}\label{void}
\end{figure}
\end{Exem}

\begin{Exem}[Exemple plus complexe] 
Faisons appel à un mot réflexif, comme \ma{lui-même}.  L'exemple est en anglais pour être moins ambigu, et pour éviter le pronom clitique \ma{se}, plus compliqué. 

\begin{tabular}{ll}
\hline
\textbf{word} & \textbf{syntactic type $u$} \\
\hline
some & $(S/(np\backslash S))/n$ \\
\hline
statements & $n$ \\
\hline
speak$\_$about & $(np\backslash S)/np$ \\
\hline
themselves & $((np\backslash S)/np)\backslash (np\backslash S)$\\
\hline
\end{tabular}

\begin{itemize} 
\item  
\ma{some} prend d'abord un nom commun puis une phrase à laquelle il manque son sujet. 
\item 
\ma{themselves} prend un verbe transitif et lui donne la catégorie des verbes intransitifs. 
\end{itemize} 

La phrase \ma{some statements speak\_about themselves} (voir dérivation figure \ref{sosta}) illustre le lien avec la logique :  le prédicat \ma{soi-même} transforme un prédicat naturellement à deux places $P(x ,y)$ en un prédicat $P(x,x)$, et donc à une seule place.  

\begin{figure}[htb]
	\begin{center}
$$
\infer[/e]{So,Sta,SpA,Refl \vdash S}{ 
\infer[/e]{So,Sta \vdash (S/(np\backslash S))}{So \vdash (S/(np\backslash S))/n & Sta \vdash n } 
& a
 \infer[\backslash e]{SpA,Refl \vdash (np \backslash S)}{ SpA \vdash (np \backslash S)/np & Refl \vdash ((np \backslash S)/np)\backslash  (np \backslash S)}}
$$	
	\end{center}
 \caption{$So$ correspond à $some$, $Sta$ à $statements$, $spA$ à $speak about$ et $Refl$ à $themselves$. Le groupe nominal  \ma{some statements} attend à sa droite une phrase privée de son sujet. Le réflexif \ma{themselves} prend à sa gauche un verbe transitif, qui demande un sujet et un objet et le transforme en verbe intransitif, ne nécessitant plus qu'un sujet. Le tout donne, par conséquent, une phrase correcte. }\label{sosta}
\end{figure}
\end{Exem}

\subsection{Sémantique de Montague}

Pour obtenir les formules représentant le sens des phrases, 
nous les allons décrire comme des $\lambda$-termes.

En plus des catégories syntaxiques,  le lexique dotera chaque mot de $\lambda$-termes 
dont les types seront une traduction des catégories syntaxiques.  
Ensuite, grâce à l'analyse syntaxique et au lexique avec ces termes ``sémantiques", nous obtiendrons un $\lambda$-terme représentant  une formule logique.

La mise en oeuvre du calcul de la représentation sémantique s'appuie sur le principe de compositionnalité cher à Frege, principe selon lequel le sens du tout est fonction du sens des parties, et de la structure syntaxique, ce  dernier ingrédient ayant été apporté par Montague --- effectivement \ma{Pierre aime Marie}  n'a pas le même sens que \ma{Marie aime Pierre}.  
Assemblons  les mots en groupes de mots , de proche en proche,  en suivant la structure syntaxique. 

Il nous faut ici préciser un peu ce qu'est un lambda terme, et comment un lambda terme peut  
représenter une formule logique. 

Les $\lambda$-termes ici considérés ont des types ; pour représenter les formules, nous utiliserons deux types de base particulier, $e$ et $t$, qui correspondent respectivement aux entités (aux individus de la logique) et aux propositions (ou aux valeurs de vérité, lorsque celles-ci sont interprétées). Ensuite, en composant les termes et les différents types, il est possible de fabriquer des fonctions et de les appliquer à des arguments.  

\begin{center}
$types ::= \eee \mid \ttt \mid types \rightarrow types$
\end{center}

\begin{itemize}
\item Pour chaque type $U$ on se donne une infinité dénombrable de variables et de constantes de type $U$ ;
\item une variable ou une constante $x$ de type $U$ est un terme de type $U$ ;
\item si $t$ est un terme de type $U\rightarrow V$\footnote{Le type $U \rightarrow V$ est le type des fonctions allant de $U$ dans $V$} et $u$ un terme de type $U$ alors $(t(u))$ est un terme de type $V$ ;
\item si $x$ est une variable de type $U$ et $t$ un terme de type $V$ alors$\lambda x.t$ est un terme de type $U \rightarrow V$ [La fonction qui à $x$ de type $U$ associe $t$ de type $V$ qui a priori dépend de $x$.].
\end{itemize}

\subsubsection{Que calcule le $\lambda$-calcul?}
Pourquoi dire $\lambda$ calcul? Pour l'instant, nous ne calculons rien avec. 
En fait, le $\lambda$-calcul permet comme son nom l'indique, de calculer, mais il le fait 
d'une manière assez stupide, qui consiste à remplacer la variable de la fonction par sa valeur.

\begin{Exem}[Calcul]
 $(\lambda x. x+2)(3+y) \rightarrow (3+y)+2.$  

Appliquer la fonction $(\lambda x. x+2)$ à $(3+y)$, c'est remplacer $x$ par $(3+y)$.
\end{Exem}

Le $\lambda$-calcul peut ``tout" calculer --- le lambda calcul non typé est l'un des modèles de la calculabilité, et le lambda calcul typé  permet de calculer ``beaucoup" de fonctions récursives totales. L'opération de base consiste à remplacer une variable par l'argument auquel la fonction est appliquée, eet cette opération est appelée  $\beta$-réduction. 
Cette réduction, dans un cadre typé, se passe bien: le calcul finit toujours, et quelque soit l'ordre dans le quel s'opèrent les réductions,  le résultat est toujours le même: tant mieux! 
  
 \subsubsection{Sémantique de Montague : syntaxe et sémantique}
 
Nous retrouverons dans ce paragraphe, l'antique (au sens propre) correspondance  qu'il y a entre une proposition et une phrase, entre un groupe nominal (défini) et un individu et d'autres du même acabit.  (voir figure \ref{corresp}). 

\begin{figure}[htb]
\begin{center}
 \begin{tabular}{|rcp{11cm}|}
 \hline
(Type syntaxique)$^*$ & = & Type sémantique \\
\hline
$S^*$ & = & $\ttt$ $\quad$ une phrase est une proposition \\
$np^*$& = & $\eee$ $\quad$ un groupe nominal est une entité\\
$n^*$& = & $\eee \rightarrow \ttt$ $\quad$ un nom commun est une partie des entités, une propriété. \\
&& Notons que noms communs et verbes intransitifs \\ 
&& se traitent de la même manière  \footnotemark.\\
$(A\backslash B)^* = (B / A)^*$ & = & $A \rightarrow B$ $\quad$ s'étend à toutes les catégories d'une grammaire \\
&&catégorielle. Quelque chose qui attendait un $A$ à sa \\
&& gauche ou à sa droite pour former un $B$ devient une fonction de $A$ dans $B$.\\
\hline 
 \end{tabular}
 \caption{On traduit les catégories en propositions, entités ou association des deux.}\label{corresp}
\end{center}
\end{figure} 

 \begin{Exem}[Verbe transitif]
Pour le verbe intransitif \ma{dormir},  le prédicat \ma{dormir} est appliqué à l'entité \ma{Pierre}. Grammaticalement, il est bien utile de savoir que l'on dit \ma{Pierre dort} et non pas \ma{dort Pierre}, mais pour la forme logique que l'on fabrique, il suffit de savoir que le prédicat \ma{dormir} s'applique à \ma{Pierre}.
 \end{Exem}
 
Ici le mot ``sujet" va être aussi bien employé au sens logique de ``sujet d'un prédicat", 
qu'au sens  grammatical de sujet d'une phrase. Ce n'est pas un hasard, cela provient de la tradition médiévale : le sujet est l'individu, le reste, c'est-à-dire le groupe verbal, est le prédicat.

 \subsubsection{Des constantes pour les opérations logiques}
 Pour exprimer la logique, nous utiliserons des constantes du lambda calcul, constantes sans lesquelles  
 il n'est point possible d'écrire de formule (voir figure \ref{constt}).

\footnotetext{Cette idée vient d'Alexandrie, environs au 3e siècle, où ont été étudiées ce genre de  correspondances. Les adjectifs posaient problème: ils partageaient des propriétés avec les verbes mais aussi avec les noms.}

\begin{figure}[htb]
\begin{center}
\begin{tabular}{r|l}
Constante & Type \\
\hline \\
$\exists$ & $(\eee\rightarrow  \ttt )\rightarrow \ttt $\\
$\forall$ & $(\eee\rightarrow \ttt)\rightarrow \ttt $\\
$\wedge$ & $\ttt\rightarrow(\ttt\rightarrow \ttt)$\\ 
$\vee$ &$\ttt\rightarrow(\ttt\rightarrow \ttt)$ \\
$\supset$ & $\ttt\rightarrow(\ttt\rightarrow \ttt)$
\end{tabular}
\caption{Traduction entre les constantes de la logique et les types en $\lambda$-calcul. Rappelons que  $e \rightarrow t$ signifie que de le lambda terme prend une entité en paramètre et renvoie une proposition. }\label{constt}
\end{center}
\end{figure}
 
 \begin{Rema}[Connecteurs et quantificateurs]
\begin{description}
\item[ $\wedge, \vee, \supset : \ttt\rightarrow (\ttt \rightarrow \ttt)$] Un connecteur logique, appliqué à une première proposition, puis à une seconde, retourne une nouvelle proposition. 
\item[$\exists, \forall: (\eee\rightarrow \ttt)\rightarrow \ttt$] L'argument d'un quantificateur, de type $(\eee \rightarrow \ttt)$ représente une propriété, un prédicat à une place --- par exemple \ma{dort}. Le quantificateur fabrique alors une proposition qui peut avoir  une valeur de vérité. Le $\forall$ prend pour argument une propriété et  fabrique la proposition qui dit si la propriété  est vraie pour tout le monde ou non (ce qui dépend de l'interprétation de la dite propriété), tandis que le $\exists$ fabrique aussi une proposition à partir d'une propriété, mais une proposition qui est vrai lorsque ladite propriété vaut d'au moins un individu. 
\end{description} 
 \end{Rema}

 \paragraph{Des constantes pour les prédicats du langage}

Essayons d'écrire ensemble quelques  formules comme des $\lambda$-termes (voir figure \ref{d51})
 \begin{Exem} Voici les $\lambda$-termes sémantiques associés à quelques mots:

 \begin{center}
 \begin{tabular}{|c|c|c|}
 \hline
 aime & $\lambda x\lambda y (aime\quad y) x$ & $x{:}\eee, y{:}\eee,aime{:}\eee \rightarrow(\eee \rightarrow \ttt)$ \\
 \hline
 \multicolumn{3}{|c|}{ « aime » \texttt{est un prédicat binaire}} \\
  \hline
 \hline
 Garance & $\lambda P (P\quad Garance)$ & $P{:}\eee\rightarrow \ttt, Garance : \eee$ \\
 \hline
 \multicolumn{3}{|c|}{ « Garance » \texttt{est décrite comme les propriétés de} « Garance » } \\
 \hline
 \end{tabular}
 \end{center}

Le terme associé au verbe $aimer$ a une particularité : il inverse les arguments, car l'usage logique est d'écrire $aime(Pierre,Marie)$ et non $aime(Marie,Pierre)$ pour signifier que \ma{Pierre aime Marie}. Il faut bien que le lambda terme associé à aimer inverse les arguments, car cette formule est construite en suivant les règles grammaticales, et celles-ci postulent que le verbe reçoit d'abord son complément d'objet avant de recevoir son sujet, le lien entre le verbe et son objet étant supposé plus intime. 
Selon une tradition issue d'Aristote, qu'on peut qualifier de ``méditerranéenne" puisqu'elle s'est ensuite propagée dans l'empire romain et byzantin notamment à Alexandrie, 
 on regroupe le verbe et ses compléments qui forme un prédicat lequel est appliqué au sujet. 
Par exemple, dans une phrase simple comme \ma{Pierre regarde Marie},
on regroupe  \ma{regarde Marie} et \ma{Pierre} est le sujet de tout le groupe verbal. 
Un prédicat aimer à deux arguments symétriques est logiquement 
possible, mais il va contre cette tradition  grammaticale, dont la linguistique moderne ne s'éloigne guère, du moins sur ce point. 
La forme sémantique de $Garance$, qui identifie un individu avec ses propriétés est à rapprocher de l'égalité de Leibniz dont nous avons déjà parlé. 
\label{d51}
 \end{Exem}
 
 \paragraph{Calcul de la forme logique : recette}
 
Etant donné l'analyse syntaxique, c'est-à-dire une preuve dans le calcul de Lambek, et un lexique contenant en sus des catégories syntaxiques les $\lambda$-termes sémantiques, la recette est assez simple: 

 \begin{enumerate}
\item Insérer sur les feuilles de l'arbre syntaxique,  c'est-à-dire les mots,  les $\lambda$-termes fournis par le lexique, en respectant les applications.
\item Réduire le $\lambda$-terme comme un groupe nominal (via $\beta$-réduction) : c'est la forme logique de l'énoncé analysé.
\end{enumerate} 
 Tous les  $\lambda$-termes obtenus qui sont de type "proposition" sont des formules, et donc transformables en types.
 
 \begin{Exem}[analyse sémantique  de \ma{some statements speak about themselves}]
Comme annoncé au début de cette partie de l'exposé,  enrichissons le tableau contenant la catégorie syntaxique $u$ uniquement en rajoutant,  le type sémantique $u^*$, ainsi que $\lambda$-terme associé $t:u^*$. 
 
\begin{center}
\begin{tabular}{|ll|}
\hline
\textbf{word} & \textbf{syntactic type $u$} \\
 &  \textbf{semantic type $u^*$} \\
 &  \textbf{semantics : $\lambda$-term of type $u^*$}  \\
 & \textbf{$x^v$ means that the variable or constant $x$ is of type $v$} \\ 
\hline
some & $(S/(np\backslash S))/n$ \\
& $(\eee  \rightarrow \ttt )\rightarrow ((\eee  \rightarrow \ttt )\rightarrow \ttt ) $ \\
&$ \lambda P^{\eee \rightarrow \ttt } \lambda Q^{\eee \rightarrow \ttt } (\exists^{(\eee \rightarrow \ttt )\rightarrow \ttt } (\lambda x^\eee (\bigwedge^{\ttt \rightarrow (\ttt \rightarrow \ttt )}(P x)(Q x))))$ \\
\hline
statements & $n$ \\
& $\eee \rightarrow \ttt $ \\
& $\lambda x^\eee (statement^{\eee \rightarrow \ttt } x)$\\
\hline
speak\_about & $(np\backslash S)/np$ \\
& $\eee  \rightarrow  (\eee  \rightarrow  \ttt )$ \\
& $\lambda y^\eee  \lambda x^\eee  ((speak\_about^{\eee \rightarrow (\eee \rightarrow \ttt )} x)y)$ \\
\hline
themselves & $((np\backslash S)/np)\backslash (np\backslash S)$\\
& $(\eee  \rightarrow (\eee  \rightarrow \ttt ))\rightarrow (\eee  \rightarrow \ttt )$ \\
& $\lambda P^{\eee \rightarrow (\eee \rightarrow \ttt )} \lambda x^\eee  ((P x)x)$\\ \hline
\end{tabular}
\end{center}

Remarquons que le nom commun et la phrase à qui il manque le sujet ont le même type : \ma{chat} et \ma{dort} sont de même type, ce sont des prédicats qui s'appliquent à une entité.

\begin{description}
\item[Some] prend deux prédicats, celui du nom commun puis celui du verbe et forme une proposition qui dit si une entité satisfait les deux prédicats à la fois. Cela se voit en analysant le $\lambda$-terme : \ma{some} appliqué à \ma{cat} et à \ma{sleep} dit s'il existe des entités \ma{qui sont des chats} et \ma{qui dorment}. 
 \item[Speak\_About] prend deux entités et retourne une valeur de vérité. Le $\lambda$-terme, inverse effectivement $x$ et $y$, en premier lieu vient de quoi on parle, est  ensuite est précisé qui en parle. 
\item[Themselves] : prend un verbe transitif et fabrique un verbe intransitif. Il prend donc une propriété qui a naturellement deux arguments, et il fabrique une proposition dont le sujet et l'objet sont identiques  (\ma{Pierre se regarde}: le pronom  \ma{se} dit que le sujet et l'objet de \ma{regarde} sont identiques). Le $\lambda$-terme a pour argument  un prédicat de type $\eee\rightarrow \eee \rightarrow \ttt$ et forme un terme de type  $(\eee \rightarrow \ttt)$, soit l'équivalent d'un prédicat à une seule place.
\end{description}

En convertissant comme indiqué les catégories syntaxiques en types sémantiques, l'analyse syntaxique présentée figure \ref{sosta}
fournit un $\lambda$-terme de ce genre : \texttt{((some statements) (themselves speak$\_$about))} de type $S^* = \ttt$.
Remplaçons  les mots par les $\lambda$-termes sémantiques correspondant et effectuons la $\beta$-réduction: 

\begin{center}
$\big((\lambda P^{\eee \rightarrow \ttt} \lambda Q^{\eee \rightarrow \ttt} ( \exists^{(\eee \rightarrow \ttt)\rightarrow \ttt } ( \lambda x^\eee (\wedge(P x)(Q x))))( \lambda x^\eee (statement^{\eee \rightarrow \ttt} x))\big)$\\[3pt] $\big(  (\lambda P^{\eee\rightarrow(\eee\rightarrow \ttt)} \lambda x^\eee ((P x)x)) ( \lambda y^\eee \lambda x^\eee ((speak\_about^{\eee \rightarrow(\eee \rightarrow \ttt)} x)y)) \big) $\\
$\downarrow \beta $\\
\vspace{3pt}
$\big(\lambda Q^{\eee\rightarrow \ttt} (\exists^{(\eee\rightarrow \ttt)\rightarrow \ttt} (\lambda x^\eee (\wedge^{\ttt\rightarrow(\ttt \rightarrow \ttt)} (statement^{\eee\rightarrow \ttt} x)(Q x)))) \big) \big( x^\eee ((speak\_about^{\eee\rightarrow(\eee \rightarrow \ttt)} x)x)\big) $\\
\vspace{3pt}
$\downarrow \beta $\\
\vspace{3pt}
$\big(\exists^{(\eee\rightarrow \ttt) \rightarrow \ttt} ( \lambda x^\eee (\wedge(statement^{\eee \rightarrow \ttt} x) ((speak\_about^{\eee \rightarrow(\eee \rightarrow \ttt)} x)x))) \big) $\\
\end{center}

Le dernier $\lambda$-terme s'écrit plus simplement ainsi:

$$\exists x{:}\eee (statement(x) \wedge speak\_about(x,x))$$

Ce qui est bien la formule logique escomptée. 
\end{Exem}

\paragraph{Résultat}

Cette interprétation  du calcul de la représentation sémantique est automatisable, et même automatisée. Elle est par exemple  implémentée dans le logiciel Grail de Richard Moot. C'est donc faisable, mais  avec de l'huile de coude, puisqu'il il faut être assez patient pour saisir manuellement toutes les entrées sémantiques du lexique (les catégories syntaxiques peuvent être acquises automatiquement). 

Dans le processus ci dessus, il n'y a aucun moyen  pour relier les sens entre eux. Par exemple, \ma{livre} est un prédicat à un argument, \ma{lire} un prédicat à deux argument, mais rien de dit qu'un \ma{livre} puisse être \ma{lu}. Rien n'est prévu  pour adapter le sens au contexte, comme \ma{classe} qui peut vouloir désigner un ensemble d'élèves mais aussi une salle dédié à l'enseignement.  Les formules logiques produites jusqu'ici ne sont pas assez subtiles, et leur production s'avère fastidieuse.
Nous allons donc prendre en compte les sortes des objets et structurer les mots en une ontologie de la langue pour faire mieux.

\section{Comment adapter le sens des mots au contexte ?}
C'est un travail relativement personnel, qui participe à une tendance actuelle où  le sens défini et calculé  précédemment est raffiné en se plaçant dans un cadre d'inspiration cognitive. 
Ces questions sont difficiles  car elles portent sur le rapport complexe qu'entretiennent la langue 
et le monde (ou l'univers du discours).  Une première constatation est qu'il faut sans doute s'éloigner de l'univers monosorte de Frege et considérer plusieurs sortes d'objets, et se diriger ainsi vers un univers plus structuré. 

La polysémie est étudiée depuis ``toujours" mais elle l'est plus particulièrement depuis les années quatre-vingts. Nous ne parlons pas ici de la polysémie de \ma{l'avocat} (fruit/profession), mais plutôt celle de \ma{livre} qui désigne tout à la fois un objet matériel et un contenu informationnel, avec un lien entre ces deux sens.  En pareil cas, il y a différents sens en rapport dont un ou plusieurs sont sélectionnés par  le ou les prédicats, et le locuteur peut vouloir affirmer des propriétés de différentes  facettes d'un même objet. Lorsque quelqu'un dit qu'un livre est  intéressant, c'est pour son contenu, mais s'il le dit volumineux  c'est en tant qu'objet matériel, mais s'il le qualifie lourd, il peut s'agir de l'une ou l'autre des deux facettes. 

\subsection{Facettes d'un mot}
Une ville peut désigner un club de foot (\ma{Marseille a battu Bordeaux}), mais aussi l'autorité qui y siège (\ma{Washington envoie des troupes en Afghanistan}), et bien sûr la ville en elle-même (\ma{Washington borde le Potomac})... De même  le mot  \ma{saumon} désigne-t-il aussi bien un animal qu'un aliment. Il est possible  d'affirmer simultanément des propriétés qui portent sur plusieurs facettes --- terme dû à Jacques Jayez --- même si ce n'est pas toujours compréhensible.

\begin{Exem}[Association de différentes facettes d'un objet]
\hspace{1em}
\begin{itemize}
\item Un livre lourd mais intéressant : aucun problème de compréhension.
\item Le diner était délicieux mais a duré des heures : idem. 
\item Un saumon rapide mais délicieux : étrange. 
\item Marseille est un grand port et bat parfois Bordeaux : incompréhensible.
\item Washington borde le Potomac et envoie des troupes en Afghanistan : utilisé uniquement de manière humoristique.
\end{itemize}
\end{Exem}

Nous verrons ci-après comment d'autres phénomènes comme le \emph{voyageur virtuel} (fictive motion) 
peuvent être traitées dans ce même cadre logique. Il s'agit de phrases comme \ma{e chemin monte pendant deux heures} qui peuvent se dire sans que personne n'emprunte ledit chemin. 

\subsubsection{Sortes et types} 

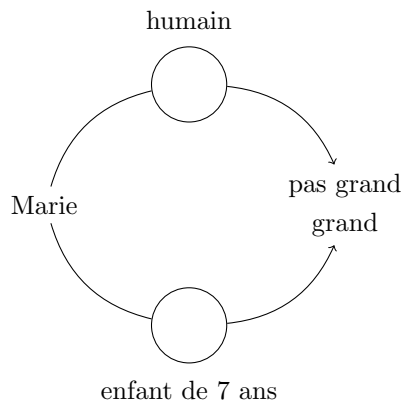
\begin{figure}[htbp]
\begin{center}
\begin{tikzpicture}[rond/.style={circle,minimum width=1cm,draw}]
	\node (m) {Marie};
	\node[above right=of m,rond] (h) {}; 
	\node[below right=of m,rond] (e) {};
	\node[below right=1cm of h,minimum width=2cm] (pg) {pas grand};
	\node[above right=1cm of e,minimum width=2cm] (g) {grand};
	\node[above=.1cm of h] {humain};
	\node[below=.1cm of e] {enfant de 7 ans};
	\path[] (m) edge[bend left] (h)
            (m) edge[bend right] (e)
            (h) edge[bend left,->] (pg) 
            (e) edge[bend right,->] (g) ;
\end{tikzpicture}
\caption{Photo de famille : Marie (7ans) est devant parce qu'elle est petite (par rapport à des adultes). Mais elle est aussi grande (pour une enfant de 7ans).}\label{F6}
\end{center}
\end{figure}

Pour exprimer ces phénomènes de facettes et de glissement de sens en contexte, 
nous utilisons une logique multisorte, $TY_n$. Cela permet de traiter les phénomènes mentionnés, tout en évitant un paradoxe qui n'existe pas dans la langue mais uniquement dans la logique monosorte de Frege et ses successeurs (voir figure \ref{F6}), paradoxe selon lequel un objet pourrait avoir des propriétés opposées. Lors d'une photo de famille, quelqu'un peut dire \ma{Marie est petite, il faut qu'elle se mette devant}, étant donné que Marie est une enfant, mais la même personne peut aussi dire \ma{Marie est grande}, et l'interlocuteur comprend alors que pour son âge, Marie est grande. 
Une analyse de ce paradoxe consiste à dire que les propriétés ne sont pas affirmées d'un objet dans l'absolu 
mais d'un objet en ``en tant que" membre d'une classe, et le paradoxe est évité 
--- de telles idées étaient déjà présentes chez Avicenne et ses élèves, notamment chez Abu-l Barak\=at al Baghd\=ad\=\i. 

Nous aurons donc plusieurs sortes d'objets, et nous nous autoriserons à quantifier sur les types.
Plusieurs systèmes permettent de telles constructions, et nous avons choisi le système F, imprédicatif mais relativement simple ---- il est assez expressif, il peut coder toute l'arithmétique intuitionniste du second ordre. 

\subsubsection{Types du second ordre}

Ce système de types inclut des variables de type et permet de quantifier sur ces variables de types. 
Par rapport aux types du système précédent, au lieu d'un type $\eee$ pour tous les individus 
nous aurons plusieurs types $\eee_i$, ce qui constitue une sorte d'ontologie plate (événements, objets physiques, contenus informationnels, êtres humains,...). 

\begin{itemize}
\item Types de base : constantes $\eee_i$ et $\ttt$, ainsi que des variables de type $\alpha \in P$.
\item Quand $T$ est un type et $\alpha$ une variable de type 
présente ou non dans $T$, $\Lambda \alpha. T$ est un type.
\item Quand $T1$ et $T2$ sont des types, $T1 \rightarrow T2$ est aussi un type.
\end{itemize}

Ces types sont les propositions quantifiées, écrite avec le seul connecteur ``$\rightarrow$".

\subsubsection{Termes du second ordre}
Les termes correspondent aux preuves du calcul des propositions quantifiées.

Cela va permettre, en présence de plusieurs types d'objets, d'appliquer une opération à tous ces différents objets. Par exemple, plutôt que d'avoir un quantificateur par type (``tous les" pour les objets physiques, ``tous les" pour les contenus informationnels, etc.), nous utiliserons  un quantificateur sur tous les types,  quantificateur qui sera ensuite  spécialisé en un quantificateur sur un type particulier. Clea permet aussi de faire des conjonctions entre des propriétés, sans même connaitre les types auxquels elles s'appliquent. 

Le système F (Girard 1971) --- dont les règles de calcul étendent naturellement la $\beta$ réduction--- est décrit par les règles si-dessous. 
$t{:}U$ (aussi noté $t^U$) signifie \ma{$t$ est un terme de type $U$}.
\begin{itemize}
\item Une variable ou une constante de type $T$,  $x : T$ ou  $x^T$ est un terme de type $T$. On dispose d'une infinité de variables de chaque type. 
\item $(f u)^V$ est un terme de type $V$ quand $u : U $ et $f : U \rightarrow V$. On appelle $(f u)$ l'application de la fonction $f$ à un terme $\tau$. 
\item $\lambda x^T. u$ est un terme de type $T \rightarrow U$ quand $x : T$ et $u : U$.
\item $t\{U\}$ est un terme de type $T[U/\alpha]$ (le type $T$ dans lequel la variable de type $\alpha$ est remplacée par $U$) quand $t : \Lambda \alpha. T$ et $U$ est un type. Un terme qui s'applique à tous les types est spécialisé pour un type $U$ particulier. Nous en verrons ci-après quelques applications. 
\item $\Lambda \alpha. t$ est un terme de type $\Lambda \alpha.T$ quand $\alpha$ est une variable de type et $t : T$ est un terme sans variable libre de type $\alpha$. Cette règle de formation construit un terme qui  généralise un terme qui fonctionne pour la variable de type $\alpha$ qui n'a rien de particulier --- non libre en hypothèse. 
\end{itemize}

Imaginons qu'il y ait des objets de type $A$, $B$, $C$ etc. 
Au lieu d'avoir un quantificateur pour les objets de type $A$, un pour ceux de type $B$ et un pour ceux de type $C$, 
c'est-à-dire des constantes $\forall^A$ de type $(A\fl \ttt)\fl \ttt$, une autre  
$\forall^B$ de type $(B\fl \ttt)\fl \ttt$ et une troisième $\forall^C$ de type $(C\fl \ttt)\fl \ttt$ 
il est possible de n'avoir qu'un seul quantificateur général  qui peut être spécialisé  à chacun des types. 

\begin{figure}[h!]
\begin{center}

\includegraphics[scale=2]{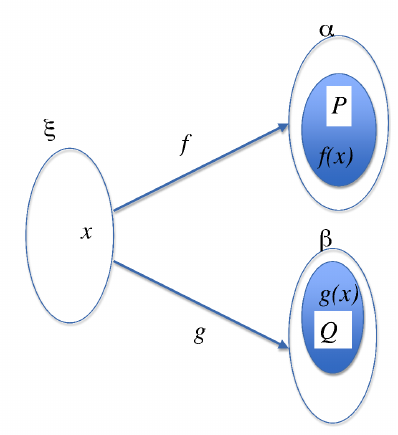} 
     \caption{Livre ($\xi$) a un sens  matériel ($\alpha$) et un sens informationnel ($\beta$). Il peut y avoir des propriétés vraies de chacune de ces deux facettes. Un livre peut par exemple être dit volumineux  ($P$) (matériellement) et intéressant  ($Q$) (d'un point de vue informationnel). Les deux informations, volumineux  et intéressant, ne s'appliquent pas à la même facette du mot \ma{livre},  mais leur conjonction  peut néanmoins être formée. 
Cette conjonction pourra être utilisée entre deux facettes de \ma{livre}, mais aussi pour lier une ville en tant que lieu géographique et en tant que représentation des habitants, etc. }\label{F7}
  \end{center}
  \end{figure}

Le même genre d'argument permet d'avoir une conjonction polymorphe qui peut coordonner deux propriétés s'appliquant à deux facettes différentes d'un même objet  (voir figure \ref{F7}). 

\begin{Exem}[La conjonction polymorphe]
Etant donnés deux prédicats $P\alpha \rightarrow t$ et $Q \beta  \rightarrow \ttt$ sur des sortes respectives d'objets $\alpha$ et $ \beta$ 
quand on a deux modifications de $ \xi $ vers $\alpha$ et vers $ \beta$ on peut coordonner des objets de type  $\xi$  :  

$\Lambda  \xi  \lambda x \xi  \lambda f \xi  \rightarrow a \lambda g \xi  \rightarrow b.(and (P (f x))(Q (g x)))$

Ceci s'applique en fait à tous les prédicats $P,Q$ et à tous  les types $\alpha, \beta$  auxquels ils s'appliquent :

$ \Lambda \alpha \Lambda  \beta  \lambda P^{\alpha \rightarrow \ttt} \lambda Q^{\beta  \rightarrow \ttt} \Lambda  \xi  \lambda x^\xi  \lambda f^{\xi  \rightarrow \alpha} \lambda g^{\xi  \rightarrow  \beta} . (and (P (f x))(Q (g x)))$
\end{Exem}
  

\subsubsection{Organisation du lexique} 

Nous ne pouvons détailler l'organisation du lexique qui permet d'automatiser le calcul des représentations sémantiques en présence de glissements de sens en contexte et qui permet ou non d'affirmer simultanément deux propriétés de différentes facettes d'un même mot. Le fait que le prédicat ne s'applique pas à  l'argument proposé correspond à une erreur de typage: $P^{a \rightarrow b} u^c$ avec $c\neq a$. Les entrées des mots correspondant à $P$ et à $u$ contiennent des $\lambda$-termes optionnels, des fonctions,  qui permettent de résoudre le conflit de type lorsque cela est souhaitable, comme pour \ma{finir un livre} où \ma{finir} attend une action. 
Le lexique spécifie si ces transformations de types sont compatibles ou au contraire s'excluent mutuellement. 

\subsubsection{Du chemin au voyageur fictif}

Dans le cadre d'un projet régional appelé Itipy, un phénomène a plus particulièrement retenu notre attention: le voyageur fictif (fictive motion). Ce phénomène apparenté aux  glissements de sens crée un lien sémantique 
entre les chemins et ceux qui les parcourent.

Le corpus du projet Itipy 
est composé de récits du 19\up{e} siècle relatant des voyages au travers des Pyrénées.   
Il contient des phrases comme \emph{... cette route monte jusqu'à Luz où l'on arrive par une jolie avenue de peupliers.}, \emph{... cette route qui monte sans cesse pendant deux heures}. 
Pour une illustration moins laborieuse de ce travail, nous utiliserons  l'exemple \emph{Le chemin monte}, non extrait du corpus mais plus rapide à traiter. La description correspond à la fois à une route et au point de vue de celui qui suit ou pourrait suivre le chemin. En effet, lorsque quelqu'un dit \ma{ce chemin monte pendant deux heures}, l'interlocuteur est obligé de penser à la personne qui suit le chemin: cela ne peut se résumer à l'altitude est une fonction croissante de l'abscisse curviligne. 

L'énoncé \ma{le chemin monte pendant deux heures}, 
pose problème car une propriété qui s'applique usuellement à un être humain, 
ou tout du moins quelque chose qui bouge, est ici appliquée à une voie. 
L'idée générale est de différencier la personne  de la voie qu'elle suit. Une transformation du terme sémantique associé au chemin, réalisée par un $\lambda$-terme optionnel va faire apparaître un $\lambda$-terme qui introduira  la personne qui suit le chemin -- quantifiée universellement : \ma{toute personne qui suit le chemin...}. 

\begin{center}
$(P^{humain\rightarrow \ttt} (u^{voie})) \qquad humain \neq voie$
\end{center}
 Clairement, c'est $u^{voie}$ qui produit un $x^{humain}$, mais si $u$
restait argument de $P$ devenu $humain$,  deux problèmes se poseraient : 
 \begin{enumerate}
\item D'une part, le quantificateur correspondant au voyageur fictif ne pourrait avoir la portée sur le prédicat.
\item D'autre part les propriétés du chemin, comme par exemple, goudronné, deviendraient des propriétés du voyageur fictif!\footnote{Ce qui est rarement le cas, hormis chez Morris et Goscinny !} La \textit{route agréable} ne pose pas de problème : il y a aussi une variable d'événement.
\end{enumerate}

En assemblant ces termes, la voie introduit un être humain. Cependant, cette opération ne porte pas sur les types atomiques, mais sur des types plus complexes c'est-à-dire sur les propriétés. 

Suivant les idées de Leibniz dont nous avons déjà parlé, 
quelque chose peut toujours être vu comme  l'ensemble des propriétés de ce quelque chose. Pour ce faire, notre modélisation utilise la montée de type (type raising): une voie peut être transformée en ensemble de propriétés des voies et c'est cet ensemble de propriétés qui vont être appliquées à l'objet de type $humain$. Le chemin ne se transforme pas en être humain, c'est le chemin ``élevé"\footnote{On appelle ça un type élevé, en référence à la montée de type : cela correspond aux propriétés qui sont vraies pour un objet.} (type raised) qui se transforme en ensemble de propriétés d'un être humain (le type humain ``élevé"). Pour des raisons techniques, il y a une variable représentant l'événement, et la  transformation s'effectue non avec $t$ (les propositions) mais avec $\tttt := \vvv \rightarrow \ttt$ (les propositions qui dépendent d'un événement, le type des événements étant $\vvv$). Le traitement est présenté de manière relativement complète en figure \ref{fictif}. Ensuite, un modificateur comme \ma{pendant deux heures} s'applique au résultat de manière habituelle.


\begin{figure} \label{fictif}
\begin{description}
\item[le:] $\Lambda \alpha \lambda P^{\alpha \rightarrow  \ttt} 
(\tau^{(\alpha \rightarrow  \ttt)\rightarrow \alpha}\ P)$ 
--- pour toute propriété des objets de type $\alpha$, \emph{le} en choisit un qui la satisfait, déterminé par le contexte.

\item[chemin:]
$\lambda x^{voie} chemin(x)$

\item[(le chemin):]

$((\Lambda \alpha \lambda P^{\alpha \rightarrow  \ttt}  (\tau^{(\alpha \rightarrow  \ttt)\rightarrow \alpha}\ P))\{voie\}\ \lambda x^{voie} chemin(x)))  $

$\quad =_{\beta} (\lambda P^{voie \rightarrow  \ttt} (\tau^{(voie\rightarrow  \ttt)\rightarrow voie}\ P)) \lambda x^{voie} chemin(x))) $

$\quad =_{\beta}  (\tau\ \lambda x^{voie} chemin(x)): voie  $

$\quad \Rightarrow  \lambda P^{voie\rightarrow \event \rightarrow  \ttt} \lambda e^{\event} (P\ (\tau\ \lambda x^{voie} chemin(x))\ e)$
\hfill  (\textit{montée de type})

\item[h:] 
$\lambda Q^{(voie\rightarrow \tttt)\rightarrow \tttt} \lambda P^{hum \rightarrow \tttt}\newline
\hspace*{1.0em} (Q\ (\lambda c^{voie}\lambda e^{\event} 
\hspace*{1.0em}  \forall(\lambda v^{hum}\ suivre(e,v,c) \imp ((P\ v) \ e))))$ 
\hfill (\textit{coercion de type})

\item[(h (le chemin)):] 

$((\lambda Q^{(voie\rightarrow \tttt)\rightarrow \tttt} 
\lambda P^{hum \rightarrow \tttt} \newline
\hspace*{1.0em} (Q\ (\lambda c^{voie}\lambda e^{\event} 
\hspace*{1.0em} 
 \forall(\lambda v^{hum }\ suivre(e,v,c) \imp ((P\ v) \ e)))))\newline
\hspace*{1.0em} (\lambda P^{voie\rightarrow \tttt} \lambda e^{\event} (P\ (\tau\ \lambda x^{voie} chemin(x))\ e)))$

\medskip 

$=_{\beta} \lambda P^{hum \rightarrow \tttt} \lambda e^{\event} 
 \forall (\lambda y^{hum } suivre(e,y,(\tau\ \lambda x^{voie} chemin(x))) \imp ((P\ x)\ e))$

\item[monte:] 
$\lambda x^{hum } \lambda e^{\event} monte(e,x)$

\item[((h (le chemin)) monte):]

$((\lambda P^{hum \rightarrow \tttt} \lambda e^{\event}\newline
\hspace*{1.0em} \forall (\lambda y^{hum } suivre(e,y,(\tau\ \lambda x^{voie} chemin(x))) \imp ((P\ x)\ e)))  (\lambda x^{hum } \lambda e^{\event} monte(e,x)))$

$\quad =_{\beta} \lambda e^{\event}\forall (\lambda y^{hum }
 suivre(e,y,(\tau\ \lambda x^{voie} chemin(x))) \imp monte(e,y))$ 
\end{description} 

\caption{Transformation d'une ``voie" (type élevé) en ``humain" (type élevé)} 
\end{figure}

\subsubsection{Comptage, quantification et individuation}
Comment compter les individus selon les différentes facettes?  Si vous avez trois exemplaires de \emph{Madame Bovary}, deux de l'\emph{\'Education sentimentale} et un exemplaire des \emph{Trois contes} vous pouvez dire \emph{J'ai descendu les livres à la cave parce que je les ai tous lus} : a priori, vous n'avez pas lu les trois copies de Madame Bovary!\footnote{On remarquera que \ma{livre} est légèrement différent de \ma{book}:
hormis pour les différents livres de la Bible, personne de dira que le livre contenant les \emph{Trois  contes} contient trois livres, tandis qu'un anglophone peut compter ce volume comme trois livres.}
Les livres peuvent être transformés en objets physiques, que quelqu'un peut compter et descendre à la cave, tandis que l'idée de lecture est liée au contenu informationnel. Des sortes de projections permettent de traiter correctement ces phénomènes. 

\subsubsection{Système F}
Pourquoi utilisons nous le système F ? Cette logique peut mettre mal à l'aise, de par sa nature imprédicative. Cela signifie que dans un type $T = \Lambda \alpha. \dots$, on peut remplacer la variable de type  $\alpha$ par $T$ lui-même,
par exemple $T = \Lambda \alpha. (\alpha \rightarrow \alpha) \rightarrow(\alpha \rightarrow \alpha)$,  peut s'appliquer à  $\alpha=T$ !

Le système F  est cependant cohérent,  en vertu  du théorème de normalisation ou de l'existence de modèles dénotationnels, et surtout, 
il bien plus simple que ses concurrents d'un point de vue syntaxique. 
Pour la tâche que nous lui assignons, à savoir la  représentation des formules de la logique multisorte,   il est très pratique: un seul 
quantificateur pour tout type (un existentiel, un universel, etc), qui peuvent être spécialisés à des types simple (objets physiques) ou plus complexe (les propositions, les verbes transitifs etc.), une conjonction polymorphe qui factorise les différentes formes de conjonction, etc. 

Cependant  n'y a pas de notion simple de sous-typage qui soit compatible avec la quantification sur les types, 
et il serait fort utile de savoir exactement ce dont a  besoin la sémantique formelle --- 
cela suggérera peut-être une notion pertinente de sous-typage. 

\subsection{Dictionnaire et univers du discours --- retour au Moyen-Âge!}

\textsl{In fine}, nos questions restent proches de celles que pouvaient se poser Aristote, Porphyre, Avicenne, Abélard, Dun Scott,... et bien d'autres. 
Quel est le lien entre l'ontologie du monde et l'ontologie de la langue ?

Un locuteur francophone peut dire \emph{ma  voiture est crevée} et non \emph{ma voiture est bouchée}, pourtant il y a des choses qui se crèvent dans une voiture, mais il y a aussi des choses qui peuvent se boucher, comme les injecteurs. Il n'empêche qu'en français, \emph{ma voiture est crevée} passe sans problème, tandis que \emph{ma voiture est bouchée} ne va pas.

Ces questions dépendent de la langue et notre modèle doit en tenir compte. 
Par exemple un anglophone ne peut dire 
\emph{my car is flat} pour\emph{a tyre of my car is flat}.

Les relations que la langue autorise,  comme la référence à certains sens particuliers, comme parler de \emph{ma voiture} pour renvoyer à l'une de ses roues, 
sont  des figures de style très présentes.  La langue ne permet pas toutes les connexions sémantiques que le monde suggère : une voiture a beau avoir des injecteurs, personé ne dira que la voiture est bouchée pour dire qu'un de ses injecteurs est bouché.

C'est une question qui peut être débattue : certains pensent que les glissements de sens sont plutôt culturels que lexicaux et que  si des personnes parlent  très longtemps des injecteurs d'une voiture, et qu'ensuite quelqu'un dit que la voiture est bouchée, il y a de fortes chances pour que l'interlocuteur comprenne qu'on parle des injecteurs. Je ne sais pas trop.

\section*{Conclusion}
Je suis désolé d'avoir un peu débordé, bien que j'aie sauté quelques paragraphes vers la fin de cette leçon. 

Ce dont j'espère vous avoir convaincu, et c'est pour ça que j'ai commencé par là, c'est que les questions posées dès l'antiquité et le moyen-âge sont très pertinentes, et toujours d'actualité.  En effet certaines n'ont pas encore de réponse satisfaisante, y compris des questions bien concrètes qui seraient fort utiles au traitement automatique des langues. 

Du côté de la syntaxe, le calcul de Lambek est très joli mais un peu restreint. Il y a alors deux solutions.
Certains, comme Moortgat gardent la structure logique de l'analyse, qui est une preuve avec les catégories, et cela conduit aux grammaires multimodales avec des postulats en plus des règles de déduction. 
Sinon, il y a la famille des modèles à deux étages : une structure profonde explique comment la phrase est structurée logiquement dont on dérive une structure de surface. Dans cette famille là, mentionnons les \textit{lambda grammars} (Muskens, 1996), les grammaires minimalistes catégorielles (Lecomte Retoré 1999), grammaires catégorielles abstraites (De Groote, 2001) grammaires de dépendances catégorielles (Dikovsky, 2001), etc.

Du point de vue de la sémantique, dont nous avons beaucoup parlé, il semble y avoir de belles avancées, en particulier du coté de la théorie des types et 
de la logique catégorique, comme l'ont montré les workshops qui se sont tenus à Oxford et à Bordeaux en novembre 2010. 

Une question intrigante est de combiner les formules multisortes du calcul des prédicats, possiblement d'ordre supérieur, avec les preuves de la  théorie des types qui assemble lesdites formules. 
Cette piste, en sémantique, est actuellement très prometteuse et très étudiée. 
Ce nouveau cadre logique, doté de nombreux types de bases et de quantifications sur les types, ouvre de nouvelles pistes sur la quantification généralisée, ainsi que sur des phénomènes difficiles comme les pluriels ou les noms massifs. 

L'interprétation dans les mondes possibles, qui semble bien établie, n'est peut-être pas aussi satisfaisante qu'il n'y paraît à première vue. 
Par exemple, les noms propres ne devraient pas être interprétés comme des constantes, afin que personne ne soit obligé  de croire que Tullius est Cicéron, mais si ce ne sont pas des constantes, encore faut-il trouver un nouveau concept plus complexe. Plus généralement, il est peu vraisemblable que nous interprétions mentalement en construisant une infinité de mondes possibles donc chacun est lui-même infini. Des interprétations des phrases  en termes de jeux et d'interaction sont sans doute plus réalistes. 

Ces deux questions ne sont que deux exemples parmi d'autres, et le domaine de la linguistique formelle a encore un bel avenir devant lui.

\providecommand{\url}[1]{\texttt{#1}}
\providecommand{\urlprefix}{URL }

\section*{Remarques bibliographiques} 

Un exposé contient peu ou pas de références, mais il convient d'en ajouter quelques unes à  la version publiée afin que le lecteur qui ne peut m'en demander puisse néanmoins en trouver et aller plus loin. 

Concernant l'histoire de la logique, la référence classique est \cite{KK86}. Sur les origines aristotéliciennes de la logique et en particulier sur les catégories qui ont traversé avec succès plus de vingt siècles nous recommandons un article relativement récent: \cite{Barnes2005cat}. 
La logique médiévale, est à l'honneur dans un ouvrage aussi captivant qu'un roman historique: \cite{libera1993philosophie}. 
Signalons aussi une bande dessinée fort bien documentée qui retrace la vie de Bertrand Russell et met en scène les logiciens du début du XXe siècle: \cite{logicomix}.
Il semble qu'elle soit un peu sévère avec Gottlob Frege, voir la critique \cite{logicomixReview}.

Concernant les grammaires catégorielles, et les extension mentionnées en conclusion,  il est tentant de faire un peu de publicité pour notre ouvrage, qui après que des versions diverses ont circulé est enfin complet et paru \cite{MootRetore2012lcg}.
L'article original de Lambek est remarquablement écrit et la lecture en est toujours recommandée \cite{Lam58}.

Concernant la sémantique de Montague mentionnons son article au titre provocateur \cite{montague:formal} ainsi que le livre collectif hollandais classique \cite{Gamut91}.
Le troisième chapitre de notre ouvrage reprend tout cela à la manière d'un cours  \cite{MootRetore2012lcg}.

La plupart des sujets à la frontière entre logique et linguistique sont fort bien présentés dans les chapitres indépendants d'environ quatre-vingt pages du \emph{Handbook of Logic and Language} qui a récemment été réédité avec des compléments très intéressants. \cite{vBtM2010}. 

Pour le lambda calcul typé, qu'il soit du premier ou du second ordre (système F) nous recommandons \cite{Girard2006pa1}.

Sur les extensions et la variation du sens des mots en contexte un ouvrage de Asher est sorti il y a peu \cite{Asher2011wow} et pour notre approche dont nous parlons à la fin de cet exposé, mentionnons \cite{BMRjolli,LMRS2012taln,Retore2012rlv}.

\bibliographystyle{plain} 

\bibliography{bigbiblio} 

\begin{thebibliography}{10}

\bibitem{Asher2011wow}
Nicholas Asher.
\newblock {\em Lexical Meaning in context -- a web of words}.
\newblock Cambridge University press, 2011.

\bibitem{Barnes2005cat}
Jonathan Barnes.
\newblock Les {C}at{\'e}gories et les \emph{{C}at{\'e}gories}.
\newblock In Otto Bruun and Lorenzo Corti, editors, {\em Les cat{\'e}gories et
  leur histoire}, Histoire de la philosophie, pages 11--80. Vrin, 2005.

\bibitem{BMRjolli}
{C}hristian {B}assac, {B}runo {M}ery, and {C}hristian {R}etor{\'e}.
\newblock {T}owards a {T}ype-{T}heoretical {A}ccount of {L}exical {S}emantics.
\newblock {\em {J}ournal of {L}ogic {L}anguage and {I}nformation},
  19(2):229--245, April 2010.
\newblock \url{http://hal.inria.fr/inria-00408308/}.

\bibitem{libera1993philosophie}
Alain de~Libera.
\newblock {\em La philosophie m{\'e}di{\'e}vale}.
\newblock Presses universitaires de France, 1993.

\bibitem{logicomix}
Ap{\'o}stolos~K. Doxi{\`a}dis, Christos Papadimitriou, Alecos Papadatos, and
  Annie~Di Donna.
\newblock {\em Logicomix}.
\newblock Vuibert, 2010.
\newblock Titre original: Logicomix: an Epic Search for Truth.

\bibitem{Gamut91}
L.~T.~F. Gamut.
\newblock {\em Logic, Language and Meaning -- Volume 2: Intensional logic and
  logical grammar}.
\newblock The University of Chicago Press, 1991.

\bibitem{Girard2006pa1}
Jean-Yves Girard.
\newblock {\em {Le point aveugle. Cours de logique. 1: Vers la perfection. 1.}}
\newblock {Paris: Hermann, \'Editeurs des Sciences et des Arts. xvi, 280~p. },
  2006.

\bibitem{KK86}
William Kneale and Martha Kneale.
\newblock {\em The development of logic}.
\newblock Oxford University Press, 3$^{\textrm{rd}}$ edition, 1986.

\bibitem{Lam58}
Joachim {L}ambek.
\newblock The mathematics of sentence structure.
\newblock {\em American mathematical monthly}, pages 154--170, 1958.

\bibitem{LMRS2012taln}
Ana{\"\i}s Lefeuvre, Richard Moot, Christian Retor{\'e}, and No{\'e}mie-Fleur
  Sandillon-Rezer.
\newblock Traitement automatique sur corpus de r{\'e}cits de voyages
  pyr{\'e}n{\'e}ens : Une analyse syntaxique, s{\'e}mantique et temporelle.
\newblock In {\em Traitement Automatique du Langage Naturel, TALN'2012},
  volume~2, pages 43--56, 2012.

\bibitem{logicomixReview}
Paolo Moncosu.
\newblock Book review: Logicomix by {A}postolos {D}oxiadis, {C}hristos {H}.
  {P}apadimitriou, {A}lecos {P}apadatos, and {A}nnie di {D}onna.
\newblock {\em Journal of Humanistic Mathematics}, 1(1):137--152, 2011.

\bibitem{montague:formal}
Richard Montague.
\newblock English as a formal language.
\newblock In Bruno Visentini, editor, {\em Linguaggi nella Societa e nella
  Tecnica}, pages 189--224. Edizioni di Communit{\`a}, Milan, Italy, 1970.
\newblock Reprinted in \cite{Mon74}.

\bibitem{MootRetore2012lcg}
Richard Moot and Christian Retor{\'e}.
\newblock {\em The logic of categorial grammars: a deductive account of natural
  language syntax and semantics}, volume 6850 of {\em LNCS}.
\newblock Springer, 2012.
\newblock
  \url{http://www.springer.com/computer/theoretical+computer+science/book/978-%
3-642-31554-1}.

\bibitem{Retore2012rlv}
Christian Retor{\'e}.
\newblock Variable types for meaning assembly: a logical syntax for generic
  noun phrases introduced by "most".
\newblock {\em Recherches Linguistiques de Vincennes}, 41:83--102, 2012.
\newblock \url{http://hal.archives-ouvertes.fr/hal-00677312}.

\bibitem{Mon74}
Richmond Thomason, editor.
\newblock {\em The collected papers of Richard Montague}.
\newblock Yale University Press, 1974.

\bibitem{vBtM2010}
Johan van Benthem and Alice ter Meulen.
\newblock {\em Handbook of Logic and Language}.
\newblock Elsevier insights. Elsevier Science, 2nd edition, 2010.

\end{thebibliography}

 \end{document}